\input ./style/arxiv-general.cfg
\documentclass[MSNbibl,number,citesort,seceqn,dvips]{arxbj}
\makeatletter
   \@ifpackageloaded{graphicx}{}{\usepackage{graphicx}}
\makeatother

\volume{22}
\issue{2}
\pubyear{2016}
\firstpage{927}
\lastpage{968}
\doi{10.3150/14-BEJ682}
\docsubty{FLA}

\makeatletter
\newcommand{\rrvert}{\vert}
\newcommand{\rrVert}{\Vert}
\newcommand{\llvert}{\vert}
\newcommand{\llVert}{\Vert}
\newcommand{\N}{\mathbb{N}}
\newcommand{\R}{\mathbb{R}}
\newcommand{\Rbar}{\overline{\R}}
\newcommand{\Cb}{\mathbb{C}}
\newcommand{\CC}{\mathcal{C}}
\newcommand{\DD}{\mathcal{D}}
\newcommand{\EE}{\mathcal{E}}
\newcommand{\FF}{\mathcal{F}}
\newcommand{\PP}{\mathcal{P}}
\newcommand{\D}{\mathbb{D}}
\newcommand{\B}{\mathbb{B}}
\newcommand{\Z}{\mathbb{Z}}
\newcommand{\Ex}{\mathrm{E}}
\newcommand{\Var}{\operatorname{Var}}
\newcommand{\Cov}{\operatorname{Cov}}
\newcommand{\MSE}{\mathrm{MSE}}
\newcommand{\IMSE}{\mathrm{IMSE}}
\newcommand{\1}{\mathbf{1}}
\renewcommand{\Pr}{\mathrm{P}}
\newtheorem{prop}{Proposition}[section]
\newtheorem{teo}[prop]{Theorem}
\newtheorem{cor}[prop]{Corollary}
\newtheorem{lem}[prop]{Lemma}
\newtheorem{cond}[prop]{Condition}
\newremark{remark}[prop]{Remark}
\makeatother

\begin{document}
\begin{frontmatter}

\title{A dependent multiplier bootstrap for the sequential empirical
copula process under strong mixing}
\runtitle{A dependent multiplier bootstrap for the sequential
empirical copula process}

\begin{aug}
\author[A]{\inits{A.}\fnms{Axel}~\snm{B\"ucher}\thanksref{A}\ead[label=e1]{axel.buecher@ruhr-uni-bochum.de}}
\and
\author[B]{\inits{I.}\fnms{Ivan}~\snm{Kojadinovic}\corref{}\thanksref{B}\ead[label=e2]{ivan.kojadinovic@univ-pau.fr}}
\address[A]{Fakult\"at f\"ur Mathematik, Ruhr-Universit\"at Bochum,
44780 Bochum, Germany.\\ \printead{e1}}
\address[B]{Laboratoire de math\'ematiques et de leurs applications,
Universit\'e de Pau et des Pays de l'Adour, UMR CNRS 5142, B.P. 1155,
64013 Pau Cedex, France. \printead{e2}}
\end{aug}

%
\received{\smonth{6} \syear{2013}}
%
\revised{\smonth{10} \syear{2014}}

%
\begin{abstract}
Two key ingredients to carry out inference on the copula of
multivariate observations are the empirical
copula process and an appropriate resampling scheme for the latter.
Among the existing techniques used for i.i.d. observations, the
multiplier bootstrap of R{\'e}millard and
Scaillet (\textit{J.~Multivariate Anal.} \textbf{100} (2009)
377--386) frequently appears to lead to inference procedures with the
best finite-sample properties.
B{\"u}cher and Ruppert (\textit{J. Multivariate Anal.} \textbf{116}
(2013) 208--229)
recently proposed an extension of this technique to strictly stationary
strongly mixing observations by adapting the \textit{dependent} multiplier
bootstrap of B\"uhlmann
(The blockwise bootstrap in time series and empirical processes (1993)
ETH Z\"urich, Section~3.3) to the empirical copula process. The main
contribution of this work is a generalization of the multiplier
resampling scheme proposed by B{\"u}cher and Ruppert along two
directions. First, the resampling scheme is now genuinely sequential,
thereby allowing to transpose to the strongly mixing setting many of
the existing multiplier tests on the unknown copula, including
nonparametric tests for change-point detection. Second, the resampling
scheme is now fully automatic as a data-adaptive procedure is proposed
which can be used to estimate the bandwidth parameter. A simulation
study is used to investigate the finite-sample performance of the
resampling scheme and provides suggestions on how to choose several
additional parameters. As by-products of this work, the validity of a
sequential version of the dependent multiplier bootstrap for empirical
processes of B\"uhlmann is obtained under weaker conditions on the
strong mixing coefficients and the multipliers, and the weak
convergence of the sequential empirical copula process is established
under many serial dependence conditions.
\end{abstract}

%
\begin{keyword}
\kwd{lag window estimator}
\kwd{multiplier central limit theorem}
\kwd{multivariate observations}
\kwd{partial-sum process}
\kwd{ranks}
\kwd{serial dependence}
\end{keyword}
\end{frontmatter}


\section{Introduction}

Let $\mathbf X$ be a $d$-dimensional random vector with continuous marginal
cumulative distribution functions (c.d.f.s) $F_1,\ldots,F_d$. From the
work of Sklar \cite{Skl59}, the c.d.f. $F$ of $\mathbf X$ can be written in a
unique way as
\[
F(\mathbf x) = C \bigl\{ F_1(x_1),\ldots,F_d(x_d)
\bigr\}, \qquad\mathbf x \in\R^d, %
\]
where the function $C\dvtx [0,1]^d \to[0,1]$ is a copula and can be
regarded as capturing the dependence among the components of $\mathbf X$.
The above equation is at the origin of the increasing use of copulas
for modeling multivariate distributions with continuous margins in many
areas such as quantitative risk management
(McNeil, Frey and Embrechts \cite{McNFreEmb05}), econometric modeling
(Patton \cite{Pat12}),
environmental modeling (Salvadori, De Michele and Kottegoda \cite
{SalDeMKotRos07}), to name a very few.

Assume that $C$ and $F_1,\ldots,F_d$ are unknown and let $\mathbf X_1,\ldots,\mathbf X_n$ be drawn from a strictly stationary sequence of continuous
$d$-dimensional random vectors with c.d.f. $F$. For any $i \in\{
1,\ldots,n\}$ and $j \in\{1,\ldots,d\}$, denote by $R_{ij}^{1:n}$ the
(mid-)rank of $X_{ij}$ among $X_{1j},\ldots,X_{nj}$ and let $\hat
U_{ij}^{1:n} = R_{ij}^{1:n}/n$. The\vspace*{1pt} random vectors $\hat{\mathbf
U}_i^{1:n} = (\hat U_{i1}^{1:n},\ldots,\hat U_{id}^{1:n})$, $i \in\{
1,\ldots,n\}$, are often referred to as \textit{pseudo-observations} from
the copula~$C$, and a natural nonparametric estimator of $C$ is the
\textit{empirical copula} of $\mathbf X_1,\ldots,\mathbf X_n$ (R{\"u}schendorf
\cite{Rus76}; Deheuvels \cite{Deh81}), frequently defined as the
empirical c.d.f. computed from the pseudo-observations, that is,
\[
C_{1:n}(\mathbf u) = \frac{1}{n} \sum_{i=1}^n
\1 \bigl( \hat{\mathbf U}_i^{1:n} \leq\mathbf u\bigr), \qquad\mathbf u
\in[0,1]^d. %
\]

The empirical copula plays a key role in most nonparametric inference
procedures on~$C$. Examples of its use for parametric inference,
nonparametric testing and goodness-of-fit testing can be found in Tsukahara \cite
{Tsu05}, R{\'e}millard and
Scaillet \cite{RemSca09},
Genest, R{\'e}millard and Beaudoin
\cite{GenRemBea09}, respectively, among many others. The asymptotics
of such procedures typically follow from the asymptotics of the \textit{empirical copula process}. With applications to change-point detection
in mind, a generalization of the latter process central to this work is
the \textit{two-sided sequential empirical copula process}. It is defined,
for any $(s,t) \in\Delta= \{ (s,t) \in[0,1]^2:   s \le t\}$ and $\mathbf u
\in[0,1]^d$,~by
%
\begin{equation}
\label{eqseqempcop} \Cb_n(s,t, \mathbf u) = \frac{1}{\sqrt{n}} \sum
_{i=\lfloor ns \rfloor
+1}^{\lfloor nt \rfloor} \bigl\{ \1 \bigl( \hat{\mathbf
U}_i^{\lfloor ns
\rfloor+1:\lfloor nt \rfloor} \leq u\bigr) - C(\mathbf u) \bigr\},
\end{equation}
where, for any $y \geq0$, $\lfloor y \rfloor$ is the greatest integer
smaller or equal than $y$. The latter process can be rewritten in terms
of the empirical copula $C_{\lfloor ns \rfloor+1:\lfloor nt \rfloor}$
of the sample $\mathbf X_{\lfloor ns \rfloor+1},\ldots,\mathbf X_{\lfloor nt
\rfloor}$ as
\[
\Cb_n(s, t, \mathbf u) = \sqrt{n} \lambda_n(s,t) \bigl\{
C_{\lfloor ns
\rfloor+1:\lfloor nt \rfloor}(\mathbf u) - C(\mathbf u) \bigr\}, \qquad(s, t, \mathbf u) \in\Delta
\times[0,1]^d, %
\]
where $\lambda_n(s,t) = (\lfloor nt \rfloor-\lfloor ns \rfloor) / n$
and with the convention that $C_{k:  k-1}(\mathbf u) = 0$ for all $\mathbf u \in
[0,1]^d$ and all $k\in\{1,\ldots, n\}$.

The quantity $\Cb_n(0, 1, \cdot, \cdot)$ is the standard empirical
copula process which has been extensively studied in the literature
(see, e.g., R{\"u}schendorf \cite{Rus76};
Gaenssler and Stute \cite{GanStu87};
Tsukahara \cite{Tsu05};
van~der Vaart and Wellner \cite{vanWel07};
Segers \cite{Seg12};
B{\"u}cher and Volgushev \cite{BucVol13}).
Notice that the process $\Cb_n(0,\cdot,\cdot,\cdot)$ does not
coincide with the sequential process initially studied by R{\"u}schendorf \cite{Rus76}
and defined by
%
\begin{equation}
\label{eqRusch} \Cb_n^\circ(s, \mathbf u) = \frac{1}{\sqrt{n}}
\sum_{i=1}^{\lfloor ns
\rfloor} \bigl\{ \1 \bigl( \hat{\mathbf
U}_i^{1:n} \leq\mathbf u\bigr) - C(\mathbf u) \bigr\}, \qquad(s,
\mathbf u) \in[0,1]^{d+1}.
\end{equation}
The above process, unlike $\Cb_n(0,\cdot,\cdot,\cdot)$, cannot be
rewritten in terms of the empirical copula unless $s=1$. Note that the
weak convergence of $\Cb_n^\circ$ was further studied by B{\"u}cher
and Volgushev \cite{BucVol13} under a large number of serial
dependence scenarios and under mild smoothness conditions on the copula.

As mentioned earlier, a first key ingredient of many of the existing
inference procedures on the unknown copula $C$ is the process $\Cb_n$
defined in~(\ref{eqseqempcop}). A second key ingredient is typically
some resampling scheme allowing to obtain replicates of $\Cb_n$. When
dealing with independent observations, several such resampling schemes
for the empirical copula process $\Cb_n(0, 1, \cdot, \cdot)$ were
proposed in the literature, ranging from the multinomial bootstrap of
Fermanian, Radulovi{\'c} and Wegkamp
\cite{FerRadWeg04} to the multiplier technique introduced in Scaillet \cite
{Sca05} and investigated further in R{\'e}millard and
Scaillet \cite{RemSca09}. Their finite-sample properties were compared
in B{\"u}cher and Dette \cite{BucDet10} who concluded that the multiplier bootstrap of R{\'e}millard and
Scaillet \cite{RemSca09} has, overall, the best finite-sample
behavior. In the case of strongly mixing observations, B{\"u}cher and
Ruppert \cite{BucRup13} recently proposed a similar resampling scheme
by adapting the \textit{dependent} multiplier bootstrap of B\"uhlmann
(\cite{Buh93}, Section~3.3) to the process $\Cb_n^\circ$ defined
in~(\ref{eqRusch}). Their empirical investigations indicate that the
latter outperforms in finite samples a block bootstrap based on the
work of K{\"u}nsch \cite{Kun89} and B{\"u}hlmann \cite{Buh94}.
Note that the idea of \textit{dependent multipliers} appearing in B\"uhlmann (\cite{Buh93},
Section~3.3) can also be found in Chen and Fan (\cite{CheFan99},
Section~5.1) and was recently independently rediscovered by Shao \cite
{Sha10} in the context of the \textit{smooth function model} but not in
the empirical process setting. For the sample mean as statistic of
interest, the latter author connected this resampling technique to the
\textit{tapered block bootstrap} of Paparoditis and Politis \cite{PapPol01}.

The main aim of this work is to provide an extended version of the
multiplier resampling scheme of B{\"u}cher and Ruppert \cite{BucRup13}
adapted to the two-sided sequential process $\Cb_n$ defined in~(\ref
{eqseqempcop}). The influence of the parameters of the resulting
bootstrap procedure is studied in detail, both theoretically and by
means of extensive simulations. An important contribution of the paper
is an approach for estimating the key bandwidth parameter which plays a
role somehow analogous to that of the block length in the block
bootstrap. As a practical consequence, the resulting dependent
multiplier technique for $\Cb_n$ can be used in a fully automatic way
and many of the existing multiplier tests on the unknown copula $C$
derived in the case of i.i.d. observations can be transposed to the
strongly mixing case. In addition, due to its sequential nature, the
resampling scheme can be used to derive nonparametric tests for
change-point detection particularly sensitive to changes in the copula.
This last point will be discussed in more detail in Section~\ref
{secdepmultCn}, and is also the subject of a companion paper (B{\"
u}cher \textit{et al}. \cite{BucKojRohSeg14}). Finally, the obtained results
could be used to develop statistical inference procedures for Markovian
copula time series models as introduced in Darsow, Nguyen and
Olsen \cite{DarNgu92}. Based on
recent results from Beare \cite{Bea10} on the mixing properties of these
time series, one could, for instance, apply the proposed multiplier
bootstrap to derive uniform confidence bands for the empirical copula
or to develop tests for simple goodness-of-fit hypotheses on the copula
in theses models.

There are two important by-products of this work that can be of
independent interest. First, the validity of a sequential version of
the dependent multiplier bootstrap for empirical processes of B\"
uhlmann (\cite{Buh93}, Section~3.3) (which has also been considered in
B{\"u}cher and Ruppert \cite{BucRup13}, proof of Proposition~2) is obtained under weaker
conditions on the rate of decay of the strong mixing coefficients and
the multipliers. The derived result is based on a sequential
unconditional multiplier central limit theorem for the multivariate
empirical process indexed by lower-left orthants that is adapted to the
case of strongly mixing observations. Second, the weak convergence of
the two-sided sequential empirical copula process $\Cb_n$ is
established under many serial dependence scenarios, including mild
strong mixing conditions.

The paper is organized as follows. The second section presents a
sequential extension of the seminal work of
B\"uhlmann (\cite{Buh93}, Section~3.3). In the third section, the
asymptotics of the two-sided sequential empirical copula process $\Cb
_n$ are obtained under many serial dependence conditions. Based on the
results of the second and third sections, a dependent multiplier
bootstrap for $\Cb_n$ is derived next. In the fifth section, the
practical steps necessary to carry out the derived bootstrap are
examined. In particular, a procedure for estimating the key bandwidth
parameter of the dependent multiplier bootstrap is proposed by adapting
to the empirical process setting the approach put forward in Politis
and White \cite{PolWhi04} and Patton, Politis and White \cite{PatPolWhi09}, among others. In
addition, two ways of generating dependent multiplier sequences central
to this resampling technique are discussed. The last section partially
reports the results of large-scale Monte Carlo experiments whose aim
was to investigate the influence in finite samples of the various
parameters involved in the dependent multiplier bootstrap for $\Cb_n$.

The following notation is used in the sequel. The arrow~``$\leadsto$''
denotes weak convergence in the sense of Definition~1.3.3 in van der Vaart and Wellner \cite{vanWel96}, and, given a set $T$, $\ell^\infty(T)$
(resp., $\CC(T)$) represents the space of all bounded (resp., continuous) real-valued functions on~$T$ equipped with the uniform metric.

\section{A dependent multiplier bootstrap for the multivariate empirical process under strong mixing}\label{secmult}

The multiplier bootstrap of R{\'e}millard and
Scaillet \cite{RemSca09} that has been adopted as a resampling
technique in the case of i.i.d. observations in many tests on the
unknown copula $C$ is a consequence of the multiplier central limit
theorem for empirical processes (see, e.g., Kosorok \cite{Kos08}, Theorem 10.1 and Corollary 10.3). A sequential version of the previous result can
be proved (see Holmes, Kojadinovic and Quessy \cite{HolKojQue13}, Theorem~1) by using the method of
proof adopted in
van der Vaart and Wellner (\cite{vanWel96}, Theorem 2.12.1). While
investigating the block bootstrap for empirical processes constructed
from strongly mixing observations, B\"uhlmann (\cite{Buh93},
Section~3.3) obtained what resembles to a conditional version of the
multiplier central limit theorem, subsequently also referred to as a
\textit{dependent multiplier bootstrap} (note that a sequential
version of this result appears in the proof of Proposition~2 of B{\"u}cher and Ruppert \cite{BucRup13}). The main idea of B\"uhlmann is to
replace i.i.d. multipliers by suitable serially dependent multipliers.
In the rest of the paper, we say that a sequence of random variables
$(\xi_{i,n})_{i \in\Z}$ is a \textit{dependent multiplier sequence} if:
\begin{longlist}[(M1)]
\item[(M1)] The sequence $(\xi_{i,n})_{i \in\Z}$ is
strictly stationary with $\Ex(\xi_{0,n}) = 0$, $\Ex(\xi_{0,n}^2) =
1$ and $\sup_{n \geq1} \Ex(\llvert  \xi_{0,n}\rrvert  ^\nu) < \infty$ for all
$\nu\geq1$, and is independent of the available sample $\mathbf X_1,\ldots,\mathbf X_n$.\vspace*{1pt}

\item[(M2)] There exists a sequence $\ell_n \to\infty$ of
strictly positive constants such that $\ell_n = \mathrm{o}(n)$ and the sequence
$(\xi_{i,n})_{i \in\Z}$ is $\ell_n$-dependent, that is, $\xi
_{i,n}$ is independent of $\xi_{i+h,n}$ for all $h > \ell_n$ and $i
\in\N$.

\item[(M3)]  There exists a function $\varphi\dvtx \R\to
[0,1]$, symmetric around 0, continuous at $0$, satisfying $\varphi
(0)=1$ and $\varphi(x)=0$ for all $\llvert  x\rrvert   > 1$ such that $\Ex(\xi_{0,n}
\xi_{h,n}) = \varphi(h/\ell_n)$ for all $h \in\Z$.
\end{longlist}

To state the main result of this section, we need to introduce
additional notation and definitions. Let $\mathbf U_1,\ldots,\mathbf U_n$ be
the unobservable sample obtained from $\mathbf X_1,\ldots,\mathbf X_n$ by the
probability integral transforms $U_{ij} = F_j(X_{ij})$, $i \in\{
1,\ldots,n\}$, $j \in\{1,\ldots,d\}$. It follows that $\mathbf U_1,\ldots,\mathbf U_n$ is a marginally uniform $d$-dimensional sample from the
unknown c.d.f.~$C$. The corresponding sequential empirical process is
then defined as
%
\begin{equation}
\label{eqseqep} \tilde\B_n(s, \mathbf u) = \frac{1}{\sqrt{n}} \sum
_{i=1}^{\lfloor ns
\rfloor} \bigl\{\1(\mathbf U_i \leq\mathbf u)
- C(\mathbf u) \bigr\}, \qquad(s,\mathbf u) \in [0,1]^{d+1}.
\end{equation}
Note that, in the rest of the paper, the notation of most of the
quantities that are directly computed from the unobservable sample $\mathbf
U_1,\ldots,\mathbf U_n$ will involve the symbol ``$\sim$.''

Furthermore,\vspace*{1pt} let $M$ be a large integer and let $(\xi_{i,n}^{(1)})_{i
\in\Z},\ldots,(\xi_{i,n}^{(M)})_{i \in\Z}$ be $M$ independent
copies of the same dependent multiplier sequence. Then, for any $m \in
\{1,\ldots,M\}$ and $(s, \mathbf u) \in[0,1]^{d+1}$,~let
%
\begin{equation}
\label{eqtildeBnm} \tilde\B_n^{(m)} (s, \mathbf u) =
\frac{1}{\sqrt{n}} \sum_{i=1}^{\lfloor ns \rfloor}
\xi_{i,n}^{(m)} \bigl\{\1(\mathbf U_i \leq\mathbf u) - C(
\mathbf u) \bigr\}.
\end{equation}
From the previous display, we see that the bandwidth sequence $\ell_n$
defined in assumption~\textup{(M2)} plays a role somehow analogous
to that of the \textit{block length} in the block bootstrap. Two ways of
forming the dependent multiplier sequences $(\xi_{i,n}^{(m)})_{i \in
\Z}$ will be presented in Section~\ref{secgenmult}.

Finally, for the sake of completeness, let us recall the notion of \textit{strongly mixing sequence}. For a sequence of $d$-dimensional random
vectors $(\mathbf Y_i)_{i \in\Z}$, the $\sigma$-field generated by $(\mathbf
Y_i)_{a \leq i \leq b}$, $a, b \in\Z\cup\{-\infty,+\infty\}$, is
denoted by $\FF_a^b$. The strong mixing coefficients corresponding to
the sequence $(\mathbf Y_i)_{i \in\Z}$ are then defined by $\alpha_0 =
1/2$ and
\[
\alpha_r = \sup_{p \in\Z} \sup_{A \in\FF_{-\infty}^p,B\in\FF
_{p+r}^{+\infty}}
\bigl\llvert \Pr(A \cap B) - \Pr(A) \Pr(B) \bigr\rrvert, \qquad r \in \N, r > 0.
\]
The sequence $(\mathbf Y_i)_{i \in\Z}$ is said to be \textit{strongly
mixing} if $\alpha_r \to0$ as $r \to\infty$.

The following result, inspired by B\"uhlmann (\cite{Buh93},
Section~3.3), could be regarded as an extension of the multiplier
central limit theorem to the sequential and strongly mixing setting for
empirical processes indexed by lower-left orthants. Its proof is given
in Appendix~\ref{proofuncondBuh}.

%
\begin{teo}[(Dependent multiplier central limit theorem)]\label{teomultCLT}
Assume that $\ell_n = \mathrm{O}(n^{1/2 - \varepsilon})$ for some $0 <
\varepsilon< 1/2$ and that $\mathbf U_1,\ldots,\mathbf U_n$ is drawn from a
strictly stationary sequence $(\mathbf U_i)_{i \in\Z}$ whose strong
mixing coefficients satisfy $\alpha_r = \mathrm{O}(r^{-a})$, $a > 3 + 3d/2$. Then,
\[
\bigl(\tilde\B_n,\tilde\B_n^{(1)},\ldots,\tilde
\B_n^{(M)} \bigr) \leadsto \bigl(\B_C,
\B_C^{(1)},\ldots,\B_C^{(M)} \bigr)
\]
in $\{\ell^\infty([0,1]^{d+1})\}^{M+1}$, where $\B_C$ is the weak
limit of the sequential empirical process $\tilde\B_n$ defined
in~(\ref{eqseqep}), and $\B_C^{(1)},\ldots,\B_C^{(M)}$ are
independent copies of $\B_C$.
\end{teo}

Before commenting on the result and the assumptions of the above
theorem, let us state a corollary that can be regarded as an
unconditional and sequential analogue of Theorem~3.2 of B\"uhlmann
\cite{Buh93}, and may be of interest for applications of empirical
processes outside the scope of copulas. Recall that $\mathbf X_1,\ldots,\mathbf
X_n$ is drawn from a strictly stationary sequence of continuous
\mbox{$d$-}dimensional random vectors with c.d.f. $F$ and that the margins of
$F$ are denoted by $F_1,\ldots,F_d$. Then, let
\[
\Z_n(s, \mathbf x) = \frac{1}{\sqrt{n}} \sum
_{i=1}^{\lfloor ns \rfloor
} \bigl\{\1(\mathbf X_i \leq\mathbf x)
- F(\mathbf x) \bigr\}, \qquad(s, \mathbf x) \in[0,1] \times\Rbar{}^d,
\]
be the usual sequential empirical process based on the observed
sequence $\mathbf X_1,\ldots,\mathbf X_n$ and, for any $m \in\{1,\ldots,M\}$, let
\[
\hat\Z_n^{(m)} (s, \mathbf x) = \frac{1}{\sqrt{n}} \sum
_{i=1}^{\lfloor
ns \rfloor} \xi_{i,n}^{(m)} \bigl
\{\1(\mathbf X_i \leq\mathbf x) - F_n(\mathbf x) \bigr\}, \qquad(s,
\mathbf x) \in[0,1] \times\Rbar{}^d, %
\]
where $\Rbar= [-\infty,\infty]$ and $F_n$ is the empirical c.d.f. computed from $\mathbf X_1,\ldots,\mathbf X_n$. The following corollary is then
a consequence of the fact that $\Z_n (s, \mathbf x) = \tilde\B_n \{s,
F_1(x_1),\ldots,F_d(x_d)\}$ for all $(s, \mathbf x) \in[0,1] \times\Rbar{}^d$ and that, under the conditions of Theorem~\ref{teomultCLT}, for
all $m \in\{1,\ldots,M\}$,
\[
\sup_{(s, \mathbf x) \in[0,1] \times\Rbar{}^d} \bigl\llvert \hat\Z_n^{(m)}
(s, \mathbf x) - \tilde\B_n^{(m)} \bigl\{s, F_1(x_1),\ldots,F_d(x_d)\bigr\} \bigr\rrvert \stackrel{\Pr} {\to}0,
\]
a proof of which follows from the proof of Lemma~\ref{lemasymequi} in
the supplementary material (B\"ucher and Kojadinovic \cite{BucKoj14}).

%
\begin{cor}[(Dependent multiplier bootstrap for $\Z_n$)]\label{cormultZn}
Assume that $\ell_n = \mathrm{O}(n^{1/2 - \varepsilon})$ for some $0 <
\varepsilon< 1/2$ and that $\mathbf X_1,\ldots,\mathbf X_n$ is drawn from a
strictly stationary sequence $(\mathbf X_i)_{i \in\Z}$ of continuous
$d$-dimensional random vectors whose strong mixing coefficients satisfy
$\alpha_r = \mathrm{O}(r^{-a})$, $a > 3 + 3d/2$. Then,
\[
\bigl(\Z_n,\hat\Z_n^{(1)},\ldots,\hat
\Z_n^{(M)} \bigr) \leadsto \bigl(\Z_F,
\Z_F^{(1)},\ldots,\Z_F^{(M)} \bigr)
\]
in $\{\ell^\infty([0,1] \times\Rbar{}^d)\}^{M+1}$, where $\Z_F$ is
the weak limit of $\Z_n$, and $\Z_F^{(1)},\ldots,\Z_F^{(M)}$ are
independent copies of $\Z_F$.
\end{cor}

%
\begin{remark} \label{remcond}
In the literature, the ``validity'' (or ``consistency'') of a bootstrap
procedure is often shown by establishing weak convergence of
conditional laws
(see, e.g., van~der Vaart \cite{Van98}, Chapter~23). In most
theoretical developments of this type, the necessary additional step of
approximating conditional laws by simulation from the random resampling
mechanism sufficiently many times is typically omitted (van~der Vaart
\cite{Van98}, page 329). An appropriate unconditional weak convergence
result of the form of the one established in Corollary~\ref
{cormultZn} (see also Segers \cite{Seg12}, and references therein for
other examples) already includes the repetition of the random
resampling mechanism and can be used to deduce consistency of a
bootstrap procedure in many situations of practical interest. A rather
general result in that direction is provided in Proposition~F.1 of the supplementary material (B\"ucher and
Kojadinovic \cite{BucKoj14}). As an important consequence, in many
situations of practical interest, both paradigms (conditional and
unconditional) can be used, and one can choose the approach that
appears to be easiest for the particular problem at hand. In the
empirical process setting, we tend to favor the unconditional paradigm
as the usual workhorses for empirical process theory, the (extended)
continuous mapping theorem and the functional delta method, appear to
be applicable under less restrictive conditions in an unconditional
setting (see, e.g., Kosorok \cite{Kos08}, Section~10.1.4).
\end{remark}

From a practical perspective, Corollary~\ref{cormultZn} is, for
instance, a first necessary step to transpose to the strongly mixing
setting the goodness-of-fit and nonparametric change-point tests based
on empirical c.d.f.s considered in Kojadinovic and Yan \cite{KojYan12} and Holmes, Kojadinovic and Quessy \cite
{HolKojQue13}, respectively. 

We end this section by a few comments on the assumptions of
Theorem~\ref{teomultCLT} and Corollary~\ref{cormultZn}:
\begin{itemize}
\item The requirement that $\ell_n = \mathrm{O}(n^{1/2 - \varepsilon})$ for
some $0 < \varepsilon< 1/2$ is used for proving the finite-dimensional
convergence involved in Theorem~\ref{teomultCLT}, while the condition
$\alpha_r = \mathrm{O}(r^{-a})$, $a > 3 + 3d/2$, is needed for the proof of the
asymptotic equicontinuity.
\item Theorem~3.2 of B\"uhlmann \cite{Buh93} can be regarded as a nonsequential conditional analogue of Corollary~\ref{cormultZn} with
slightly more constrained multiplier random variables. The condition on
the rate of decay of the strong mixing coefficients in that result is
$\sum_{r=0}^{\infty} (r+1)^p \alpha_r^{1/2} < \infty$ with $p =
\max\{8d+12,\lfloor2/\varepsilon \rfloor + 1\}$ and is therefore
stronger than the condition involved in Theorem~\ref{teomultCLT}.

\item The condition on the strong mixing coefficients in Theorem~\ref
{teomultCLT} and Corollary~\ref{cormultZn} is clearly satisfied if
$\mathbf X_1,\ldots,\mathbf X_n$ are i.i.d., so that the above unconditional
resampling scheme remains valid for independent observations. In the
latter case however, the Monte Carlo experiments carried out in B{\"
u}cher and Ruppert \cite{BucRup13} suggest that a simpler scheme with
i.i.d. multipliers (based, e.g., on Theorem~1
of Holmes, Kojadinovic and Quessy \cite{HolKojQue13}) will lead to better finite-sample performance. As
noted by a referee, this was to be expected since the use of a
resampling scheme designed to capture dependence for observations that
are serially independent should naturally result in an efficiency loss,
especially if the tuning parameter is estimated.
\end{itemize}

\section{Asymptotics of the sequential empirical copula process under serial dependence}

In the case of i.i.d. observations, the classical empirical copula
process turns out to be asymptotically equivalent to a linear
functional of the multivariate sequential empirical process $\tilde\B
_n$ defined in~(\ref{eqseqep}) (see Segers \cite{Seg12}, Proposition
4.3). This representation is at the heart of the multiplier bootstrap
of R{\'e}millard and
Scaillet \cite{RemSca09}. Obtaining such an asymptotic representation
for the two-sided sequential empirical copula process $\Cb_n$ defined
in~(\ref{eqseqempcop}) is therefore a preliminary step before a \textit{dependent} multiplier bootstrap for $\Cb_n$ under strong mixing can be
derived as a consequence of Theorem~\ref{teomultCLT}. The desired
result is actually a corollary of a more general result. Indeed, in
this section, the asymptotics of $\Cb_n$ are established under many
serial dependence scenarios as a consequence of the weak convergence of
the multivariate sequential empirical process $\tilde\B_n$. More
specifically, the following condition is considered.

%
\begin{cond}
\label{condweak}
The sample $\mathbf U_1,\ldots,\mathbf U_n$ is drawn from a strictly stationary
sequence $(\mathbf U_i)_{i \in\Z}$ such that $\tilde\B_n$ converges
weakly in $\ell^\infty([0,1]^{d+1})$ to a tight centered Gaussian
process $\B_C$ concentrated on
\begin{eqnarray*}
&& \bigl\{ \alpha^\star\in\CC\bigl([0,1]^{d+1}\bigr)\dvtx
\alpha^\star(s,\mathbf u) = 0\mbox{ if one of the components of $(s,\mathbf u)$
is 0 and}
\\
&&\hspace*{159pt} \alpha ^\star(s,1,\ldots,1) = 0 \mbox{ for all } s \in(0,1] \bigr\}.
\end{eqnarray*}
\end{cond}

Note that, in the case of serial independence, the above condition is
an immediate consequence of Theorem 2.12.1 of van der Vaart and
Wellner \cite{vanWel96}. As shall be discussed below, it is also met
under strong mixing.

We also consider the following smoothness condition on $C$ proposed by
Segers \cite{Seg12}. As explained by the latter author, this condition
is nonrestrictive in the sense that it is necessary for the candidate
weak limit of $\Cb_n$ to exist pointwise and have continuous sample paths.

%
\begin{cond}
\label{condpd}
For any $j \in\{1,\ldots,d\}$, the partial derivatives $\dot C_j =
\partial C/\partial u_j$ exist and are continuous on $\{ \mathbf u \in
[0,1]^d:   u_j \in(0,1) \}$.
\end{cond}

As we continue, for any $j \in\{1,\ldots,d\}$, we define $\dot C_j$ to
be zero on the set $\{ \mathbf u \in[0,1]^d:   u_j \in\{0,1\} \}$
(see also Segers \cite{Seg12}; B{\"u}cher and Volgushev \cite
{BucVol13}). It then follows that, under Condition~\ref{condpd},
$\dot C_j$ is defined on the whole of $[0,1]^d$. Also, for any $j \in\{
1,\ldots,d\}$ and any $\mathbf u \in[0,1]^d$, $\mathbf u^{(j)}$ is the vector
of $[0,1]^d$ defined by $u^{(j)}_i = u_j$ if $i = j$ and 1 otherwise.

Finally, in order to study $\Cb_n$, we need to be able to easily go
back and forth between normalized ranks and empirical quantile
functions. To this end, ties must not occur. In the case of serial
independence, it is sufficient to assume that the marginal
distributions are continuous. However, in the case of serial
dependence, continuity of the marginal distributions is \textit{not}
sufficient to guarantee the absence of ties (see, e.g., B{\"u}cher and Segers \cite{BucSeg14}, Example~4.2). This leads to a last condition.

%
\begin{cond}
\label{condnoTies}
For any $j \in\{1,\ldots, d\}$, there are no ties in the component
series $X_{1j},\ldots, X_{nj}$ with probability one.
\end{cond}

The following theorem is the main result of this section. It is proved
in Appendix~\ref{proofstuterep}.

%
\begin{teo}[(Asymptotics of the sequential empirical copula process)]\label{teostuterep}
Under Conditions~\ref{condweak},~\ref{condpd} and~\ref{condnoTies},
\[
\sup_{(s,t,\mathbf u)\in\Delta\times[0,1]^d} \bigl\llvert \Cb_n(s, t, \mathbf u) -
\tilde\Cb_n(s, t, \mathbf u) \bigr\rrvert \stackrel{\Pr} {\to}0,
\]
where
%
\begin{equation}
\label{eqtildeCn} \tilde\Cb_n(s, t, \mathbf u) = \bigl\{ \tilde
\B_n(t, \mathbf u) - \tilde\B _n(s, \mathbf u) \bigr\} - \sum
_{j=1}^d \dot C_j(\mathbf u)
\bigl\{ \tilde\B_n\bigl(t, \mathbf u^{(j)}\bigr) - \tilde
\B_n\bigl(s, \mathbf u^{(j)}\bigr) \bigr\}.
\end{equation}
Consequently, $\Cb_n \leadsto\Cb_C$ in $\ell^\infty(\Delta\times
[0,1]^d)$, where, for $(s,t, \mathbf u) \in\Delta\times[0,1]^d$,
%
\begin{equation}
\label{eqec} \Cb_C(s, t, \mathbf u) = \bigl\{ \B_C(t,
\mathbf u) - \B_C(s, \mathbf u)\bigr\} - \sum_{j=1}^d
\dot C_j(\mathbf u) \bigl\{ \B_C\bigl(t, \mathbf
u^{(j)}\bigr) - \B_C\bigl(s, \mathbf u^{(j)}\bigr)
\bigr\}.
\end{equation}
\end{teo}

The asymptotics of $\Cb_n$ under strong mixing immediately follow from
the previous theorem. The necessary tool is Theorem~1 of B\"ucher \cite
{Buc14}, which states that, if $\mathbf U_1,\ldots,\mathbf U_n$ is drawn from a
strictly stationary sequence $(\mathbf U_i)_{i \in\Z}$ whose\vspace*{1pt} strong
mixing coefficients satisfy $\alpha_r = \mathrm{O}(r^{-a})$, $a>1$, then
$\tilde\B_n \leadsto\B_C$ in $\ell^\infty([0,1]^{d+1})$. In other
words, $\mathbf U_1,\ldots,\mathbf U_n$ satisfies Condition~\ref{condweak}.

%
\begin{cor}
\label{corweakCnsm}
Assume that $\mathbf X_1,\ldots,\mathbf X_n$ is drawn from a strictly
stationary sequence $(\mathbf X_i)_{i \in\Z}$ whose strong mixing
coefficients satisfy $\alpha_r = \mathrm{O}(r^{-a})$, $a > 1$. Then, under
Conditions~\ref{condpd} and~\ref{condnoTies},
\[
\sup_{(s,t,\mathbf u)\in\Delta\times[0,1]^d} \bigl\llvert \Cb_n(s, t, \mathbf u) -
\tilde\Cb_n(s, t, \mathbf u) \bigr\rrvert \stackrel{\Pr} {\to}0,
\]
where $\tilde\Cb_n$ is defined in~(\ref{eqtildeCn}).
\end{cor}

The conditions of the above corollary are, for instance, satisfied (with
much to spare) when $\mathbf X_1,\ldots,\mathbf X_n$ is drawn from a stationary
vector ARMA process with absolutely continuous innovations (see
Mokkadem \cite{Mok88}).

\section{A dependent multiplier bootstrap for $\mathbb{C}_n$ under strong mixing}
\label{secdepmultCn}

Analogously to the approach adopted in R{\'e}millard and
Scaillet \cite{RemSca09} (see also Segers \cite{Seg12}), we shall now
combine the asymptotic representation for $\Cb_n$ stated in
Corollary~\ref{corweakCnsm} with Theorem~\ref{teomultCLT} to show
the validity of a dependent multiplier bootstrap for $\Cb_n$ under
strong mixing. The corresponding result, stated in Proposition~\ref
{propmultCn} below, can be regarded as an extension of Proposition~2
in B{\"u}cher and Ruppert \cite{BucRup13}, where a similar but
conditional result was established for the process $\Cb_n^\circ$
defined in~(\ref{eqRusch}) under stricter conditions on the mixing
rate and the multipliers.

The underlying idea is as follows: the fact that the limiting vector of
processes in Theorem~\ref{teomultCLT} has independent components
suggests regarding $\tilde\B_n^{(1)},\ldots,\tilde\B_n^{(M)}$ as
``almost'' independent copies of $\tilde\B_n$ when $n$ is large.
Unfortunately, the $\tilde\B_n^{(m)}$ cannot be computed because $C$
is unknown and the sample $\mathbf U_1,\ldots,\mathbf U_n$ is unobservable.
Estimating $C$ by the empirical copula $C_{1:n}$ and $\mathbf U_1,\ldots,\mathbf U_n$ by the pseudo-observations $\hat{\mathbf U}_1^{1:n},\ldots,\hat
{\mathbf U}_n^{1:n}$, we obtain\vspace*{1pt} the following computable version of $\tilde
\B_n^{(m)}$ defined, for any $(s,\mathbf u) \in[0,1]^{d+1}$, by
%
\begin{equation}
\label{eqhatBnm} \hat\B_n^{(m)}(s,\mathbf u) = \frac{1}{\sqrt{n}}
\sum_{i=1}^{\lfloor
ns \rfloor} \xi_{i,n}^{(m)}
\bigl\{ \1 \bigl( \hat{\mathbf U}_i^{1:n} \leq\mathbf u \bigr) -
C_{1:n}(\mathbf u) \bigr\}.
\end{equation}
Starting from the asymptotic representation of $\Cb_n$ in terms of
$\tilde\B_n$ stated in Corollary~\ref{corweakCnsm}, we see that,
to obtain ``almost'' independent copies of $\Cb_n$ for large $n$ in
the spirit of R{\'e}millard and
Scaillet \cite{RemSca09}, we additionally need to estimate the partial
derivatives $\dot C_j$, $j \in\{1,\ldots,d\}$, appearing in~(\ref
{eqec}). As we continue, we consider estimators $\dot C_{j,n}$ of
$\dot C_j$ satisfying the following condition put forward in Segers
\cite{Seg12}.

%
\begin{cond}
\label{condestpd}
There exists a constant $K > 0$ such that $\llvert  \dot C_{j,n}(\mathbf u)\rrvert   \leq
K$ for all $j \in\{1,\ldots,d\}$, $n \geq1$ and $\mathbf u \in[0,1]^d$,
and, for any $\delta\in(0,1/2)$ and $j \in\{1,\ldots,d\}$,
\[
\mathop{\sup_{\mathbf u \in[0,1]^d}}_{u_j \in[\delta,1-\delta]} \bigl\llvert \dot C_{j,n}(\mathbf u) - \dot
C_j(\mathbf u)\bigr\rrvert \stackrel{\Pr} {\to}0. %
\]
\end{cond}

Three estimators of the partial derivatives satisfying Condition~\ref
{condestpd} are discussed in Section~\ref{secestimatorspd}.

We can now define empirical processes that can be fully computed and
that, under appropriate conditions, can be regarded as ``almost''
independent copies of $\Cb_n$ for large~$n$. For any $m \in\{1,\ldots,M\}$ and $(s,t,\mathbf u) \in\Delta\times[0,1]^d$, let
%
\begin{eqnarray}\label{eqhatCbnm}
\hat\Cb_n^{(m)}(s, t, \mathbf u) &=& \bigl\{ \hat
\B_n^{(m)}(t, \mathbf u) - \hat \B_n^{(m)}(s,
\mathbf u)\bigr\}
\nonumber\\[-8pt]\\[-8pt]\nonumber
&&{}- \sum_{j=1}^d \dot
C_{j,n}(\mathbf u) \bigl\{ \hat\B _n^{(m)}\bigl(t, \mathbf
u^{(j)}\bigr) - \hat\B_n^{(m)}\bigl(s, \mathbf
u^{(j)}\bigr)\bigr\}.
\end{eqnarray}

The following proposition is a consequence of Corollary~\ref
{corweakCnsm} and Theorem~\ref{teomultCLT} and can be proved by
adapting the arguments of
Segers (\cite{Seg12}, proof of Proposition 4.3) to the current
sequential and strongly mixing setting. Its proof can be found in
Section~D of the supplementary material
(B\"ucher and Kojadinovic \cite{BucKoj14}).

%
\begin{prop}[(Dependent multiplier bootstrap for $\Cb_n$)]\label{propmultCn}
Assume that $\ell_n = \mathrm{O}(n^{1/2 - \varepsilon})$ for some $0 <
\varepsilon< 1/2$ and that $\mathbf X_1,\ldots,\mathbf X_n$ is drawn from a
strictly stationary sequence $(\mathbf X_i)_{i \in\Z}$ whose strong
mixing coefficients satisfy $\alpha_r = \mathrm{O}(r^{-a})$, $a > 3 + 3d/2$.
Then, under Conditions~\ref{condpd},~\ref{condnoTies} and~\ref{condestpd},
\[
\bigl(\Cb_n,\hat\Cb_n^{(1)},\ldots,\hat
\Cb_n^{(M)} \bigr) \leadsto \bigl(\Cb_C,
\Cb_C^{(1)},\ldots,\Cb_C^{(M)} \bigr)
\]
in $\{\ell^\infty(\Delta\times[0,1]^d)\}^{M+1}$, where $\Cb_C$ is
the weak limit of the two-sided sequential empirical copula process
$\Cb_n$ defined in~(\ref{eqec}), and $\Cb_C^{(1)},\ldots,\Cb
_C^{(M)}$ are independent copies of $\Cb_C$.
\end{prop}

We end this section by briefly illustrating how Proposition~\ref
{propmultCn} can be used in the context of change-point detection. As
discussed in B{\"u}cher \textit{et al}. \cite{BucKojRohSeg14}, a broad class of
nonparametric tests for change-point detection particularly sensitive
to changes in the copula can be derived from the process
\[
\D_n(s,\mathbf u) = \sqrt{n} \lambda_n(0,s)
\lambda_n(s,1) \bigl\{ C_{1:\lfloor ns \rfloor}(\mathbf u) - C_{\lfloor ns \rfloor+1:n}(
\mathbf u) \bigr\}, \qquad(s,\mathbf u) \in[0,1]^{d+1}. %
\]
The above definition is a mere transposition to the copula context of
the ``classical construction'' adopted, for instance, in
Cs{\"o}rg{\H{o}} and Horv{\'a}th (\cite{CsoHor97}, Section~2.6).
Under the null hypothesis of no change in the distribution, the process
$\D_n$ can be simply rewritten as
\[
\D_n(s,\mathbf u) = \lambda_n(s,1) \Cb_n(0,s,\mathbf
u) - \lambda_n(0,s) \Cb_n(s,1,\mathbf u), \qquad(s,\mathbf u)
\in[0,1]^{d+1}. %
\]
To be able to compute approximate $p$-values for statistics derived
from $\D_n$ (given the unwieldy nature of the weak limit of $\D_n$),
it is then natural to define the processes
\[
\hat\D_n^{(m)}(s,\mathbf u) = \lambda_n(s,1) \hat
\Cb_n^{(m)}(0,s,\mathbf u) - \lambda_n(0,s) \hat
\Cb_n^{(m)}(s,1,\mathbf u), \qquad(s,\mathbf u) \in
[0,1]^{d+1}, %
\]
$m \in\{1,\ldots,M\}$, which could be thought of as ``almost''
independent copies of $\D_n$ under the null hypothesis of no change in
the distribution. Under the null and the conditions of Proposition~\ref
{propmultCn}, we immediately obtain from Proposition~\ref
{propmultCn} that $\D_n, \hat\D_n^{(1)},\ldots, \hat\D_n^{(M)}$
weakly converge jointly to independent copies of the same limit. As
discussed in Remark~\ref{remcond}, the latter result is the key step
for establishing that classical tests based on $\D_n$ hold their level
asymptotically. To illustrate this point further, let us focus on the
Kolmogorov--Smirnov statistic $W_n = \sup_{(s,\mathbf u) \in[0,1]^{d+1}}
\llvert   \D_n(s,\mathbf u)\rrvert  $ and let $\hat W_n^{(m)} = \sup_{(s,\mathbf u) \in
[0,1]^{d+1}} \llvert   \hat\D_n^{(m)}(s,\mathbf u)\rrvert  $, $m \in\{1,\ldots,M\}$. The
continuous mapping\vspace*{1pt} theorem then implies that, under the null and the
conditions of Proposition~\ref{propmultCn}, $(W_n, W_n^{(1)},\ldots,
W_{n}^{(M)}) \leadsto( W, W^{(1)},\ldots, W^{(M)})$, where\vspace*{1pt} $W$, the
weak limit of $W_n$, is a continuous random variable, and
$W^{(1)},\ldots,W^{(M)}$ are independent copies of $W$. The above
unconditional result ensures that the conclusion of Proposition~F.1
in Section~F of the
supplementary material (B\"ucher and Kojadinovic \cite{BucKoj14})
holds, which implies that a test based on $W_n$ whose approximate
$p$-value is computed as $M^{-1} \sum_{m=1}^M \1 (\hat{W}_n^{(m)}
\geq W_n )$ will hold its level asymptotically as $n \to\infty$
followed by $M \to\infty$. To show that such a test is consistent
under the alternative of changes in the copula only, one typically\vspace*{1pt}
needs to prove that $n^{-1/2}W_n \stackrel{\Pr}{\to}c>0$ and that,
for any $m \in\{1,\ldots,M\}$, $W_n^{(m)} = \mathrm{O}_\Pr(\ell_n^{1/2})$,
also under the alternative (see, e.g., Inoue \cite{Ino01}, Theorem 2.5
for related results in the context of nonparametric change-point
detection in multivariate c.d.f.s).

Additional details, simulation results as well as illustrations on
financial data can be found in B{\"u}cher \textit{et al}. \cite{BucKojRohSeg14}
for tests based on maximally selected Cram\'er--von Mises statistics.




\section{Practical issues}

The practical use of the derived dependent multiplier bootstrap for
$\Cb_n$ requires the generation of dependent multiplier sequences and
the estimation of the partial derivatives of the copula. These two
practical issues are discussed in the second and third subsection
below, while the first subsection addresses the key choice of the
bandwidth parameter $\ell_n$ involved in the definition of dependent
multiplier sequences.

\subsection{Estimation of the bandwidth parameter \texorpdfstring{$\ell_n$}{$ell_n$}}\label{secband}

The bandwidth parameter $\ell_n$ defined in assumption~\textup{(M2)} plays a role somehow similar to that of the block length in
the block bootstrap of K{\"u}nsch \cite{Kun89}. Its value is therefore expected
to have a crucial influence on the finite-sample performance of the
dependent multiplier bootstrap for $\Cb_n$. The choice of a similar
bandwidth parameter is discussed, for instance, in Paparoditis and
Politis \cite{PapPol01} for the tapered block bootstrap using results
from K{\"u}nsch \cite{Kun89}. Related results are presented in B\"uhlmann
(\cite{Buh93}, Lemmas~3.12 and~3.13) and Shao (\cite{Sha10},
Proposition 2.1) for the dependent multiplier bootstrap when the
statistic of interest is the sample mean. The aim of this section is to
extend the aforementioned results to the dependent multiplier bootstrap
for $\Cb_n$ and propose an estimator of $\ell_n$ in the spirit of
those investigated in Paparoditis and Politis \cite{PapPol01}, Politis
and White \cite{PolWhi04} and Patton, Politis and White \cite{PatPolWhi09} for other resampling
schemes. Since the dependent multiplier bootstrap for $\Cb_n$ is based
on the corresponding bootstrap approximation for $\tilde\B_n$, we
propose to base our estimator of the bandwidth parameter on the
accuracy of the latter technique.

Let $\Ex_\xi$ and $\Cov_\xi$ denote the expectation and covariance,
respectively, conditional on the data $\mathbf X_1,\ldots,\mathbf X_n$, and,
for any $\mathbf u, \mathbf v \in[0,1]^d$, let $\sigma_C(\mathbf u, \mathbf v) = \Cov
\{ \B_C(1, \mathbf u), \B_C(1, \mathbf v) \}$. Now, fix $m \in\{1,\ldots,M\}
$ and, for any $\mathbf u, \mathbf v \in[0,1]^d$, let
%
\begin{eqnarray}\label{eqtildesig}
\nonumber
\tilde\sigma_n(\mathbf u, \mathbf v) &=& \Cov_\xi\bigl
\{ \tilde\B_n^{(m)}(1, \mathbf u), \tilde\B_n^{(m)}(1,
\mathbf v)\bigr\}
\\
&=& \Ex_\xi\bigl\{ \tilde\B _n^{(m)}(1,
\mathbf u) \tilde\B_n^{(m)}(1, \mathbf v)\bigr\}
\nonumber\\[-8pt]\\[-8pt]\nonumber
&=& \frac{1}{n} \sum_{i,j=1}^n
\Ex_\xi\bigl( \xi_{i,n}^{(m)}
\xi_{j,n}^{(m)} \bigr) \bigl\{\1(\mathbf U_i \leq\mathbf
u) - C(\mathbf u) \bigr\} \bigl\{\1(\mathbf U_j \leq\mathbf v) - C(\mathbf v) \bigr
\}
\\
&=& \frac{1}{n} \sum_{i,j=1}^n
\varphi\bigl\{ (i-j) / \ell_n \bigr\} \bigl\{\1(\mathbf U_i
\leq\mathbf u) - C(\mathbf u) \bigr\} \bigl\{\1(\mathbf U_j \leq\mathbf v) - C(\mathbf v)
\bigr\},\nonumber
\end{eqnarray}
where $\tilde\B_n^{(m)}$ is defined in~(\ref{eqtildeBnm}). For the
moment, although it is based on the unobservable sample $\mathbf U_1,\ldots,\mathbf U_n$ and the unknown copula $C$, we shall regard $\tilde\sigma
_n(\mathbf u, \mathbf v)$ as an estimator of $\sigma_C(\mathbf u, \mathbf v)$.

The following two results extend Lemmas~3.12 and 3.13 of B\"uhlmann
\cite{Buh93} and Proposition~2.1 of Shao \cite{Sha10}. They can be
proved by adapting the arguments used in the proofs of Lemmas~3.12 and
3.13 of B\"uhlmann \cite{Buh93}. The resulting proofs are given in the
supplementary material
(B\"ucher and Kojadinovic \cite{BucKoj14}) for completeness.

%
\begin{prop}
\label{propbias}
Assume that $\ell_n = \mathrm{O}(n^{1/2 - \varepsilon})$ for some $0 <
\varepsilon< 1/2$, that $\mathbf U_1,\ldots,\mathbf U_n$ is drawn from a
strictly stationary sequence $(\mathbf U_i)_{i \in\Z}$ whose strong
mixing coefficients satisfy $\alpha_r = \mathrm{O}(r^{-a})$, $a > 3$, and that
$\varphi$ defined in assumption~\textup{(M3)} is additionally
twice continuously differentiable on $[-1,1]$ with $\varphi''(0) \neq
0$. Then, for any $\mathbf u, \mathbf v \in[0,1]^d$,
\[
\Ex\bigl\{ \tilde\sigma_n(\mathbf u, \mathbf v) \bigr\} -
\sigma_C(\mathbf u, \mathbf v) = \frac{\Gamma(\mathbf u, \mathbf v)}{\ell_n^2} + r_{n,1}(\mathbf u,
\mathbf v), %
\]
where $\sup_{\mathbf u, \mathbf v \in[0,1]^d} \llvert   r_{n,1}(\mathbf u, \mathbf v) \rrvert   =
\mathrm{o}(\ell_n^{-2})$ and
\[
\Gamma(\mathbf u, \mathbf v) = \frac{\varphi''(0)}{2} \sum_{k=-\infty
}^\infty
k^2 \gamma(k,\mathbf u, \mathbf v)\qquad\mbox{with } \gamma(k, \mathbf u, \mathbf v)
= \Cov\bigl\{ \1(\mathbf U_0 \leq\mathbf u), \1(\mathbf U_k \leq\mathbf
v) \bigr\}. %
\]
\end{prop}

%
\begin{prop}
\label{propvar}
Assume that $\mathbf U_1,\ldots,\mathbf U_n$ is drawn from a strictly
stationary sequence $(\mathbf U_i)_{i \in\Z}$ whose strong mixing
coefficients satisfy $\alpha_r = \mathrm{O}(r^{-a})$, $a > 3$, and that there
exists $\lambda> 0$ such that $\varphi$ defined in assumption~\textup{(M3)} additionally satisfies $\llvert  \varphi(x) - \varphi(y)\rrvert   \leq
\lambda\llvert  x - y \rrvert  $ for all $x,y \in\R$. Then, for any $\mathbf u, \mathbf v \in
[0,1]^d$,
\[
\Var\bigl\{ \tilde\sigma_n(\mathbf u, \mathbf v) \bigr\} = \frac{\ell_n}{n}
\Delta (\mathbf u,\mathbf v) + r_{n,2}(\mathbf u, \mathbf v), %
\]
where
\[
\Delta(\mathbf u,\mathbf v) = \biggl\{ \int_{-1}^1
\varphi(x)^2 \,\mathrm{d}x \biggr\} \bigl[ \sigma_C(\mathbf u,
\mathbf u) \sigma_C(\mathbf v,\mathbf v) + \bigl\{ \sigma_C(\mathbf u,
\mathbf v) \bigr\}^2 \bigr] %
\]
and $\sup_{\mathbf u,\mathbf v \in[0,1]^d} \llvert   r_{n,2}(\mathbf u, \mathbf v) \rrvert   = \mathrm{o}(\ell_n/n)$.
\end{prop}

Under the combined conditions of Propositions~\ref{propbias} and~\ref
{propvar}, we have that, for any $\mathbf u, \mathbf v \in[0,1]^2$, the mean
squared error of $\tilde\sigma_n(\mathbf u, \mathbf v)$ is
\[
\MSE\bigl\{ \tilde\sigma_n(\mathbf u, \mathbf v) \bigr\} = \frac{\{ \Gamma(\mathbf u,
\mathbf v) \}^2 }{\ell_n^4}
+ \Delta(\mathbf u, \mathbf v) \frac{\ell_n}{n} + r_n(\mathbf u, \mathbf v),
\]
where $r_n(\mathbf u, \mathbf v) = \{r_{n,1}(\mathbf u, \mathbf v)\}^2 + 2 \Gamma(\mathbf
u, \mathbf v) r_{n,1}(\mathbf u, \mathbf v)/\ell_n^2 + r_{n,2}(\mathbf u, \mathbf v)$.
This allows us to define the integrated mean squared error
%
\begin{equation}
\label{eqimse} \IMSE_n = \int_{[0,1]^{2d}} \MSE\bigl\{
\tilde\sigma_n(\mathbf u, \mathbf v) \bigr\} \,\mathrm{d}\mathbf u \,\mathrm{d}\mathbf v
\sim\frac{ \bar\Gamma^2 }{\ell
_n^4} + \bar\Delta\frac{\ell_n}{n},
\end{equation}
where
%
\begin{equation}
\label{eqbarGD} \bar\Gamma^2 = \int_{[0,1]^{2d}} \bigl\{
\Gamma(\mathbf u, \mathbf v) \bigr\}^2 \,\mathrm{d}\mathbf u \,\mathrm{d}\mathbf v \quad
\mbox{and}\quad\bar\Delta= \int_{[0,1]^{2d}} \Delta(\mathbf u, \mathbf v)
\,\mathrm{d}\mathbf u \,\mathrm{d}\mathbf v.
\end{equation}
Notice that $\bar\Delta$ can be rewritten as
\[
\bar\Delta= \biggl\{ \int_{-1}^1
\varphi(x)^2 \,\mathrm{d}x \biggr\} \biggl[ \biggl\{ \int
_{[0,1]^d} \sigma_C(\mathbf u, \mathbf u) \,\mathrm{d}\mathbf u
\biggr\}^2 + \int_{[0,1]^{2d}} \bigl\{
\sigma_C(\mathbf u, \mathbf v) \bigr\} ^2 \,\mathrm{d}\mathbf u
\,\mathrm{d}\mathbf v \biggr]. %
\]
Differentiating the function $x \mapsto\bar\Gamma^2/x^4 + \bar
\Delta x/n$ and equating the derivative to zero, we obtain that the
value of $\ell_n$ that minimizes $\IMSE_n$ is, asymptotically,
%
\begin{equation}
\label{eqlnopt} \ell_n^{\mathrm{opt}} = \biggl( \frac{4 \bar\Gamma^2 }{\bar\Delta}
\biggr)^{1/5} n^{1/5}.
\end{equation}

%
\begin{figure}

\includegraphics{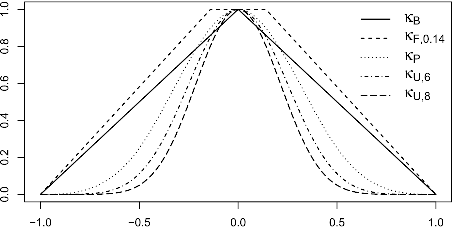}

\caption{Graphs of the functions $\kappa_B$, $\kappa_{F,0.14}$ and
$\kappa_P$, as well as $\kappa_{U,6}$ and $\kappa_{U,8}$ defined in
Section~\protect\ref{subseccov}.}\label{figkernels}
\end{figure}

From~(\ref{eqlnopt}), we see that, to estimate $\ell_n^{\mathrm{opt}}$, we
need to estimate the infinite sums $K(\mathbf u, \mathbf v)=\sum_{k \in\Z}
k^2 \gamma(k, \mathbf u, \mathbf v)$ and $\sigma_C(\mathbf u, \mathbf v) = \sum_{k
\in\Z} \gamma(k, \mathbf u, \mathbf v)$ for all $\mathbf u,\mathbf v \in[0,1]^d$.
Should $\mathbf U_1,\ldots,\mathbf U_n$ be observable, this could be done by
adapting the procedures described in Paparoditis and Politis (\cite
{PapPol01}, page 1111) or
Politis and White (\cite{PolWhi04}, Section~3) to the current
empirical process setting. Let $L \geq1$ be an integer to be
determined from $\mathbf X_1,\ldots,\mathbf X_n$ later and fix $\mathbf u,\mathbf v \in
[0,1]^d$. Proceeding in the spirit of Politis and Romano \cite{PolRom95} and Politis \cite
{Pol03}, the quantity $K(\mathbf u, \mathbf v)$ could be estimated by $\check
K_n(\mathbf u, \mathbf v) = \sum_{k=-L}^L \kappa_{F,0.5}(k/L) k^2 \check
\gamma_n(k,\mathbf u, \mathbf v)$,
where
%
\begin{equation}
\label{eqflattop} \kappa_{F,c}(x) = \bigl[ \bigl\{ \bigl(1-\llvert x
\rrvert \bigr)/(1-c) \bigr\} \vee0 \bigr] \wedge1, \qquad c \in[0,1],
\end{equation}
is the ``flat top'' (trapezoidal) kernel parametrized by $c \in[0,1]$
(see Figure~\ref{figkernels}), and $\check\gamma_n(k,\mathbf u, \mathbf v)$
is the estimated cross-covariance at lag $k \in\{-(n-1),\ldots,n-1\}$,
computed from the sequences $\{\1(\mathbf U_i \leq\mathbf u)\}_{i \in\{
1,\ldots,n\}}$ and $\{\1(\mathbf U_i \leq\mathbf v)\}_{i \in\{1,\ldots,n\}}$,
that is,
\[
\check\gamma_n(k,\mathbf u, \mathbf v) = \cases{\displaystyle n^{-1} \sum
_{i=1}^{n - k} \bigl\{\1(\mathbf U_i
\leq\mathbf u) - \tilde H_n(\mathbf u) \bigr\} \bigl\{\1(\mathbf U_{i+k}
\leq\mathbf v) - \tilde H_n(\mathbf v) \bigr\}, &\quad $k \geq0$,
\vspace*{5pt}\cr
\displaystyle n^{-1} \sum_{i=1-k}^n \bigl\{\1(
\mathbf U_i \leq\mathbf u) - \tilde H_n(\mathbf u) \bigr\} \bigl\{\1(
\mathbf U_{i+k} \leq\mathbf v) - \tilde H_n(\mathbf v) \bigr\}, &\quad
$k \leq0$,}
\]
with $\tilde H_n$ being the empirical c.d.f. computed from $\mathbf
U_1,\ldots,\mathbf U_n$. Similarly, $\sigma_C(\mathbf u, \mathbf v)$ could be
estimated by
\[
\check\sigma_n(\mathbf u, \mathbf v) = \sum_{k=-L}^L
\kappa_{F,0.5}(k/L) \check\gamma_n(k,\mathbf u, \mathbf v).
\]
As\vspace*{1pt} $\mathbf U_1,\ldots,\mathbf U_n$ is unobservable, it is natural to consider
the sample of pseudo-observations $\hat{\mathbf U}^{1:n}_1,\ldots,\hat
{\mathbf U}^{1:n}_n$ instead, and to replace $\check\gamma_n(k,\mathbf u, \mathbf
v)$ by
\[
\hat\gamma_n(k,\mathbf u, \mathbf v) = \cases{\displaystyle n^{-1} \sum
_{i=1}^{n - k} \bigl\{\1\bigl(\hat{\mathbf
U}^{1:n}_i \leq\mathbf u\bigr) - C_{1:n}(\mathbf u) \bigr
\} \bigl\{\1\bigl(\hat{\mathbf U}^{1:n}_{i+k} \leq\mathbf v\bigr) -
C_{1:n}(\mathbf v) \bigr\}, &\quad $k \geq0$,
\vspace*{5pt}\cr
\displaystyle n^{-1} \sum
_{i=1-k}^n \bigl\{\1\bigl(\hat{\mathbf
U}^{1:n}_i \leq\mathbf u\bigr) - C_{1:n}(\mathbf u) \bigr
\} \bigl\{\1\bigl(\hat{\mathbf U}^{1:n}_{i+k} \leq\mathbf v\bigr) -
C_{1:n}(\mathbf v) \bigr\}, &\quad $k \leq0$,}
\]
which gives the computable estimators
%
\begin{eqnarray}\label{eqhatsigmaKn}
\hat\sigma_n(\mathbf u, \mathbf v) &=& \sum
_{k=-L}^L \kappa_{F,0.5}(k/L) \hat
\gamma_n(k,\mathbf u, \mathbf v) \quad\mbox{and}
\nonumber\\[-8pt]\\[-8pt]\nonumber
\hat K_n(\mathbf u, \mathbf v) &=& \sum_{k=-L}^L
\kappa_{F,0.5}(k/L) k^2 \hat\gamma_n(k,\mathbf u, \mathbf v)
\end{eqnarray}
of\vspace*{2pt} $\sigma_C(\mathbf u, \mathbf v)$ and $\sum_{k \in\Z} k^2 \gamma(\mathbf u,
\mathbf v)$, respectively.

To estimate $\bar\Gamma^2$ and $\bar\Delta$ defined in~(\ref
{eqbarGD}), we then propose to use a grid $\{\mathbf u_i\}_{i \in\{
1,\ldots,g\}}$ of $g$ points uniformly spaced over $(0,1)^d$, and to compute
\[
\hat{\bar\Gamma}^2_n = \frac{\{ \varphi''(0) \}^2}{4}
\frac
{1}{g^2} \sum_{i,j=1}^g \bigl\{
\hat K_n(\mathbf u_i, \mathbf u_j) \bigr
\}^2 %
\]
and
\[
\hat{\bar\Delta}_n = \biggl\{ \int_{-1}^1
\varphi(x)^2 \,\mathrm{d}x \biggr\} \Biggl( \Biggl\{
\frac{1}{g} \sum_{i=1}^g \hat
\sigma _n(\mathbf u_i, \mathbf u_i) \Biggr
\}^2 + \frac{1}{g^2} \sum_{i,j=1}^g
\bigl\{ \hat\sigma_n(\mathbf u_i, \mathbf u_j)
\bigr\}^2 \Biggr), %
\]
respectively. Plugging\vspace*{1pt} these into~(\ref{eqlnopt}), we obtain an
estimator of $\ell_n^{\mathrm{opt}}$ which shall be denoted as $\hat\ell
_n^{\mathrm{opt}}$ as we continue.

The above estimator depends on the choice of the integer $L$ appearing
in~(\ref{eqhatsigmaKn}). To estimate~$L$, we suggest proceeding
along the lines of Politis and White (\cite{PolWhi04}, Section~3.2)
(see also Paparoditis and Politis \cite{PapPol01}, page 1112). Let
$\hat\rho_j(k)$, $j \in\{1,\ldots,d\}$, be the autocorrelation
function at lag $k$ estimated from the sample $X_{1j},\ldots,X_{nj}$.
For any $j \in\{1,\ldots,d\}$, let $L_j$ be the smallest integer after
which $\hat\rho_j(k)$ appears negligible. Notice that the latter can
be determined automatically by means of the algorithm described in
detail in Politis and White (\cite{PolWhi04}, Section~3.2). Our
implementation is based on Matlab code by A.J. Patton (available on his
web page) and its \textsf{R} version by J. Racine and C. Parmeter.
Then, we merely suggest taking $L = 2\psi(L_1,\ldots,L_d )$, where
$\psi$ is some aggregation function such as the median, the mean, the
minimum or the maximum. The previous approach is clearly not the only
possible multivariate extension of the procedure of Politis and White
\cite{PolWhi04}. 
Nonetheless, the choice $\psi= \mathrm{median}$ was found to give
meaningful results in our Monte Carlo experiments partially reported in
Section~\ref{secmc}.

\subsection{Generation of dependent multiplier sequences}\label{secgenmult}

The practical use of the results stated in Sections~\ref{secmult}
and~\ref{secdepmultCn} requires the generation of dependent
multiplier random variables satisfying assumptions~\textup{(M1)}, \textup{(M2)} and~\textup{(M3)}. We
describe two ways of constructing such dependent sequences. The first
one generalizes the moving average approach proposed by B\"uhlmann
(\cite{Buh93}, Section~6.2) (see also B{\"u}cher and Ruppert \cite
{BucRup13}) and produces multipliers that satisfy assumption~\textup{(M3)} only asymptotically. The second one was suggested by
Shao \cite{Sha10} and is based on the calculation of the square root
of the covariance matrix implicitly defined in assumption~\textup{(M3)}.

\subsubsection{The moving average approach}\label{subsecma}

Let $\kappa$ be some positive bounded real function symmetric around
zero such that $\kappa(x) > 0$ for all $\llvert  x\rrvert   < 1$. Let $b_n$ be a
sequence of integers such that $b_n \to\infty$, $b_n = \mathrm{o}(n)$ and $b_n
\geq1$ for all $n \in\N$. Let $Z_1,\ldots,Z_{n+2b_n-2}$ be i.i.d. random variables independent of the available sample $\mathbf X_1,\ldots,\mathbf X_n$ such that $\Ex(Z_1)=0$, $\Ex(Z_1^2) = 1$ and $\Ex
(\llvert  Z_1\rrvert  ^\nu) < \infty$ for all $\nu\geq1$. Then, let $\ell
_n=2b_n-1$ and, for any $j \in\{1,\ldots,\ell_n\}$, let $w_{j,n} =
\kappa\{(j-b_n)/b_n\}$ and $\tilde w_{j,n} = w_{j,n} ( \sum_{j'=1}^{\ell_n} w_{j',n}^2 )^{-1/2}$. Finally, for all $i \in\{
1,\ldots,n\}$, let
\[
\xi_{i,n} = \sum_{j=1}^{\ell_n}
\tilde w_{j,n} Z_{j+i-1}. %
\]
Clearly, $\xi_{1,n},\ldots,\xi_{n,n}$ are identically distributed
with $\Ex(\xi_{1,n}) = 0$, $\Ex(\xi_{1,n}^2) = 1$ and it can be
verified that $\sup_{n \geq1} \Ex(\llvert  \xi_{1,n}\rrvert  ^\nu) < \infty$ for
all $\nu\geq1$. 
Furthermore, $\xi_{1,n},\ldots,\xi_{n,n}$ are $(\ell_n-1)$-dependent
and, for any $i \in\{1,\ldots,n\}$ and $r \in\{0,\ldots,(\ell_n-1)
\wedge n\}$,
\begin{eqnarray*}
\Cov(\xi_{i,n} \xi_{i+r,n}) &=& \sum
_{j=1}^{\ell_n} \sum_{j'=1}^{\ell_n}
\tilde w_{j,n} \tilde w_{j',n} \Ex (Z_{j+i-1}Z_{j'+i+r-1})
= \sum_{j=r+1}^{\ell_n} \tilde w_{j,n}
\tilde w_{j-r,n}
\\
&=& \Biggl( \sum_{j=1}^{\ell_n}
w_{j,n}^2 \Biggr)^{-1} \sum
_{j=r+1}^{\ell_n} \kappa\bigl\{(j-b_n)/b_n
\bigr\} \kappa\bigl\{(j-r-b_n)/b_n\bigr\}.
\end{eqnarray*}

For practical reasons, only a sequence of size $n$ has been generated.
From the previous developments, we immediately have that the infinite
size version of $\xi_{1,n},\ldots,\xi_{n,n}$
satisfies assumptions~\textup{(M1)}~and~\textup{(M2)} (as $(\ell_n-1)$-dependence clearly implies $\ell
_n$-dependence). Let us now verify that it satisfies assumption~\textup{(M3)} asymptotically.

Assume additionally that $\kappa(x) = 0$ for all $\llvert  x\rrvert   > 1$, and, for
any $f,g\dvtx \Z\to\R$, let $f * g$ denote the discrete convolution of
$f$ and $g$, that is, $f * g(r) = \sum_{j=-\infty}^\infty f(j)
g(r-j)$, $r \in\Z$. Then, let $\kappa_{b_n}(j) = \kappa(j/b_n)$, $j
\in\Z$, and notice that the previous covariance can be written as
\[
\Cov(\xi_{i,n} \xi_{i+r,n}) = \frac{\sum_{j=-\infty}^\infty
\kappa_{b_n}(j - b_n) \kappa_{b_n}(j-r - b_n)}{\kappa_{b_n} * \kappa
_{b_n}(0)} + \mathrm{o}(1)
= \frac{ \kappa_{b_n} * \kappa_{b_n}(r) }{ \kappa_{b_n} * \kappa
_{b_n}(0) } + \mathrm{o}(1) %
\]
for all $i \in\{1,\ldots,n\}$ and $r \in\{0,\ldots,n-i\}$, where the
$\mathrm{o}(1)$ term comes from the fact that $\kappa(1)$ is not necessarily
equal to 0.

Assume furthermore that there exists $\lambda> 0$ such that $\llvert  \kappa
(x) - \kappa(y)\rrvert   \leq\lambda\llvert  x - y \rrvert  $ for all $x,y \in[-1,1]$ and
let $r_n$ be a positive sequence such that $r_n/b_n \to\gamma\in
[0,1]$. We shall now check that $b_n^{-1} \kappa_{b_n} * \kappa
_{b_n}(r_n) \to\kappa\star\kappa(\gamma)$, where $\star$ denotes
the convolution operator between real functions. We have
\[
\frac{1}{b_n} \kappa_{b_n} * \kappa_{b_n}(r_n)
= \frac{1}{b_n} \sum_{j=-b_n}^{b_n}
\kappa(j/b_n) \kappa\bigl\{(r_n-j)/b_n\bigr\}.
\]
On the one hand,
\begin{eqnarray*}
&& \Biggl\llvert \frac{1}{b_n} \sum_{j=-b_n}^{b_n}
\kappa(j/b_n) \kappa\bigl\{ (r_n-j)/b_n\bigr
\} -\frac{1}{b_n} \sum_{j=-b_n}^{b_n}
\kappa(j/b_n) \kappa(\gamma-j/b_n) \Biggr\rrvert
\\
&&\quad \leq \lambda\llvert r_n/b_n
- \gamma\rrvert \frac{2b_n+1}{b_n} \sup_{x \in\R} \kappa(x) \to0,
\end{eqnarray*}
and, and on the other hand,
\[
\frac{1}{b_n} \sum_{j=-b_n}^{b_n}
\kappa(j/b_n) \kappa(\gamma -j/b_n) \to\int
_{-1}^1 \kappa(x) \kappa(\gamma-x) \,\mathrm{d}x =
\kappa\star\kappa(\gamma). %
\]
It follows that
\[
\frac{ \kappa_{b_n} * \kappa_{b_n}(r_n) }{ \kappa_{b_n} * \kappa
_{b_n}(0) } \to\frac{ \kappa\star\kappa(\gamma) }{ \kappa\star
\kappa(0) }. %
\]
Now, let
%
\begin{equation}
\label{eqkappavarphi} \varphi(x) = \frac{\kappa\star\kappa(2x)}{\kappa\star\kappa
(0)}, \qquad x \in\R,
\end{equation}
where the factor 2 ensures that $\varphi(x) = 0$ for all $\llvert  x\rrvert   > 1$.
Then, for large $n$, $\Cov(\xi_{i,n} \xi_{j,n}) \approx\varphi\{
(i-j)/\ell_n\}$, for any $i,j \in\{1,\ldots,n\}$. Hence, the infinite
size version of $\xi_{1,n},\ldots,\xi_{n,n}$ satisfies
assumption~\textup{(M3)} asymptotically.


In our numerical experiments, we considered several popular \textit{kernels} for the function~$\kappa$ (see, e.g., Andrews \cite{And91}),
defined, for any $x \in\R$, as
\begin{eqnarray*}
\mbox{Truncated:}&\qquad& \kappa_{T}(x) = \1\bigl(\llvert x\rrvert \leq1
\bigr),
\\
\mbox{Bartlett:}&\qquad& \kappa_{B}(x) = \bigl(1 - \llvert x\rrvert
\bigr) \vee0,
\\
\mbox{Parzen:}&\qquad& \kappa_{P}(x) = \bigl(1 -
6x^2 + 6\llvert x\rrvert ^3\bigr) \1\bigl(\llvert x
\rrvert \leq 1/2\bigr) + 2\bigl(1-\llvert x\rrvert \bigr)^3 \1
\bigl(1/2 < \llvert x\rrvert \leq1\bigr),
\end{eqnarray*}
as well as the flat top kernel already defined in~(\ref{eqflattop}).
The above kernels satisfy all the assumptions on the function $\kappa$
mentioned previously. Their graphs are represented in Figure~\ref
{figkernels}. 
The flat top (or trapezoidal) kernel, parametrized by $c \in[0,1]$,
was used in Paparoditis and Politis \cite{PapPol01} in the context of
the tapered block bootstrap for the mean. These authors found that,
within the class of trapezoidal kernels symmetric around 0.5 and with
support $(0,1)$, $\kappa_{F,0.14}$, rescaled and shifted to have
support $(0,1)$, minimizes the asymptotic mean squared error of the
bootstrapping procedure. The latter kernel was also used in Shao \cite
{Sha10} who connected the tapered block bootstrap with the dependent
multiplier bootstrap for the mean.

\subsubsection{The covariance matrix approach}\label{subseccov}

Let $\ell_n$ be a sequence of strictly positive constants such that
$\ell_n \to\infty$ and $\ell_n = \mathrm{o}(n)$. Let $\varphi$ be a
function satisfying assumption~\textup{(M3)} such that,
additionally, $\int_{-\infty}^\infty\varphi(u) e^{-\mathrm{i}ux} \,\mathrm{d}u \geq0$ for all $x \in\R$, and let $\bolds\Sigma_n$ be the $n
\times n$ (covariance) matrix whose elements are defined by $\varphi\{
(i-j)/\ell_n\}$, $i,j \in\{1,\ldots,n\}$. The integral condition on
$\varphi$ ensures\vspace*{1pt} that $\bolds\Sigma_n$ is positive definite which in
turn ensures the existence of $\bolds\Sigma_n^{1/2}$. From a practical
perspective, $\bolds\Sigma_n^{1/2}$ can be computed either by
diagonalization, singular value decomposition or Cholesky factorization
of $\bolds\Sigma_n$. We use the first approach. Then, let $Z_1,\ldots,Z_n$ be i.i.d. standard normal random variables independent of the
available sample $\mathbf X_1,\ldots,\mathbf X_n$. A dependent multiplier
sequence $\xi_{1,n}, \ldots,\xi_{n,n}$ can then be simply obtained as
\[
[\xi_{1,n},\ldots, \xi_{n,n}]^\top= \bolds
\Sigma_n^{1/2} [Z_1,\ldots,Z_n]^\top.
\]
If $\varphi(1) > 0$, then the above construction generates $\ell
_n$-dependent multipliers, while if $\varphi(1) = 0$, the generated
sequence is $(\ell_n-1)$-dependent. Clearly, the infinite size version
of $\xi_{1,n},\ldots,\xi_{n,n}$ satisfies assumptions~\textup{(M1)},~\textup{(M2)} and~\textup{(M3)}.

From a practical perspective, for the function $\varphi$, we
considered the Bartlett and Parzen kernels $\kappa_B$ and $\kappa_P$,
as well as $\kappa_{U,6}$ and $\kappa_{U,8}$, where $\kappa_{U,p}$
is the density function of the sum of $p$ independent uniforms centered
at~0, normalized so that it equals 1 at~0, and rescaled to have support
$(-1,1)$. The functions $\kappa_{U,6}$ and $\kappa_{U,8}$ are
represented in Figure~\ref{figkernels}. Notice that $\kappa_T =
\kappa_{U,1}$, $\kappa_B = \kappa_{U,2}$ and $\kappa_P = \kappa
_{U,4}$. This also implies that $\kappa_{U,8}$ is a rescaled and
normalized version of the convolution of $\kappa_P$ with itself, that
is, $\kappa_{U,8}(x) = \kappa_P \star\kappa_P(2x)/\kappa_P \star
\kappa_P(0)$ for all $x \in\R$. A numerically stable and efficient
way of computing $\kappa_{U,p}$ consists of using \textit{divided differences}
(see, e.g., Agarwal, Dalpatadu and Singh \cite{AgaDalSin02}). Finally,
note that the truncated and flat top kernels cannot be used as they do
not satisfy the integral condition ensuring that $\Sigma_n$ is
positive definite.

%
\begin{remark}
\label{remrel}
In the case of the moving average approach presented in Section~\ref
{subsecma}, we have seen that $\kappa$ determines $\varphi$
asymptotically through~(\ref{eqkappavarphi}). It follows that, for
an initial standard normal i.i.d. sequence, the same value of $\ell
_n$ and for large $n$, we could expect the dependent multiplier
sequences generated by the moving average and the covariance matrix
approaches, respectively, to give close results when $\kappa$ in
Section~\ref{subsecma} and $\varphi$ in Section~\ref{subseccov}
are related through~(\ref{eqkappavarphi}). For instance, all other
parameters being similar, using the Bartlett kernel for $\kappa$ in
Section~\ref{subsecma} should produce similar results to using the
Parzen kernel for $\varphi$ in Section~\ref{subseccov}.
\end{remark}

\subsection{Estimation of the partial derivatives of the copula}\label{secestimatorspd}

For the estimators of the partial derivatives appearing in~(\ref
{eqhatCbnm}), we considered three possible definitions proposed in the
literature. The first one is that of R{\'e}millard and
Scaillet \cite{RemSca09} who suggested to estimate the partial
derivatives $\dot C_j$, $j \in\{1,\ldots,d\}$, by finite-differences as
%
\begin{eqnarray}
\label{eqpdestRS} \dot C_{j,n}(\mathbf u) &=& \frac{1}{2 n^{-1/2}} \bigl\{
C_n \bigl(u _1,\ldots,u_{j-1},u_j
+ n^{-1/2},u_{j+1},\ldots,u_d \bigr)
\nonumber\\[-8pt]\\[-8pt]\nonumber
&&\hspace*{33pt}{}-
C_n \bigl( u_1,\ldots,u_{j-1},u_j -
n^{-1/2},u_{j+1},\ldots,u_d \bigr) \bigr\}, \qquad
\mathbf u \in[0,1]^d.\qquad
\end{eqnarray}
A slightly different definition consisting of a ``boundary correction''
was proposed in Kojadinovic, Segers and Yan (\cite{KojSegYan11}, page
706). Yet another definition is mentioned in B{\"u}cher and Ruppert
(\cite{BucRup13}, page 212).\vspace*{1pt} Note that, for any $\delta\in(0,1/2)$,
all three definitions coincide on the set $\{\mathbf u \in[0,1]^d:   u_j
\in[\delta,1-\delta]\}$ provided $n$ is taken large enough. Now,
under the assumptions of Corollary~\ref{corweakCnsm}, we have that
$\Cb_n(0,1,\cdot) \leadsto\Cb_C(0,1,\cdot)$ in $\ell^\infty
([0,1]^d)$. The latter weak convergence implies the first statement of
Lemma~2 of
Kojadinovic, Segers and Yan \cite{KojSegYan11}, which in turn implies
that Condition~\ref{condestpd} is satisfied for the above defined
$\dot C_{j,n}$ as well as for the two slightly different definitions
considered in Kojadinovic, Segers and Yan (\cite{KojSegYan11},
page~706) and B{\"u}cher and Ruppert (\cite{BucRup13}, page 212), respectively.

\section{Monte Carlo experiments}\label{secmc}

To investigate the finite-sample performance of the proposed dependent
multiplier bootstrap, we considered several statistics derived from the
sequential empirical copula process $\Cb_n$ defined in~(\ref
{eqseqempcop}). With applications to statistical tests in mind, we
mostly focus in this section on Cram\'er--von-Mises and
Kolomogorov--Smirnov statistics obtained from $\Cb_n(0,1,\cdot)$.
Results for some simpler functionals can be found in Section~G of the supplementary material
(B\"ucher and Kojadinovic \cite{BucKoj14}).

Recall that $M$ is a large integer, and let
%
\begin{eqnarray}\label{eqSn}
S_n &=& \int_{[0,1]^d} \bigl\{
\Cb_n(0,1,\mathbf u) \bigr\}^2 \,\mathrm{d}\mathbf u \quad
\mbox{and}\quad
\nonumber\\[-8pt]\\[-8pt]\nonumber
S_n^{(m)} &=& \int_{[0,1]^d}
\bigl\{ \hat\Cb _n^{(m)}(0,1,\mathbf u) \bigr\}^2
\,\mathrm{d}\mathbf u, \qquad m \in\{1,\ldots,M\},
\end{eqnarray}
where $\hat\Cb_n^{(m)}$ is defined in~(\ref{eqhatCbnm}) with the
partial derivative estimators defined as discussed later in this
section. Under the\vspace*{1pt} conditions of Proposition~\ref{propmultCn} and
from the continuous mapping theorem, we then immediately have that
$(S_n,S_n^{(1)},\ldots,S_n^{(M)})$ converges weakly to
$(S,S^{(1)},\ldots,S^{(M)})$, where $S = \int_{[0,1]^d} \{ \Cb
_C(0,1,\mathbf u) \}^2 \,\mathrm{d}\mathbf u$ and $S^{(1)},\ldots,S^{(M)}$ are
independent copies of $S$.\vspace*{1pt}

The\vspace*{1pt} first aim of our Monte Carlo experiments was to assess the quality
of the estimation of the quantiles of $S$ by the empirical quantiles of
the sample $S_n^{(1)},\ldots,S_n^{(M)}$. Let $S_n^{(1:  M)} \leq\cdots
\leq S_n^{(M:  M)}$ denote the corresponding order\vspace*{1pt} statistics. An
estimator of the quantile of $S$ of order $p \in(0,1)$ is then simply
$S_n^{(\lfloor pM \rfloor:  M)}$. For each data generating scenario, the
target theoretical quantiles of $S$ of order $p$ were accurately
estimated empirically from $10^5$ realizations of $S_{1000}$ for $p \in
\PP= \{0.25,0.5,0.75,0.9,0.95,0.99\}$. Then, for\vspace*{1pt} each data generating
scenario, $N=1000$ samples $\mathbf X_1,\ldots,\mathbf X_n$ were generated and,
for each sample, $S_n^{(\lfloor pM \rfloor:  M)}$ was computed for each
$p \in\PP$ using the dependent multiplier bootstrap with $M=2500$
yielding, for each $p \in\PP$, $N$ estimates of the quantile of $S$
of order $p$. This allowed us to compute, for each data generating
scenario and each $p \in\PP$, the empirical bias and the empirical
mean squared error (MSE) of the estimators of the quantiles of $S$ of
order $p$. Similar simulations were performed for the
Kolmogorov--Smirnov statistic. Specifically, let
%
\begin{equation}
\label{eqTn} \hspace*{-18pt}
T_n = \sup_{\mathbf u \in[0,1]^d} \bigl\llvert
\Cb_n(0,1,\mathbf u)\bigr\rrvert \quad\mbox {and}\quad
T_n^{(m)} = \sup_{\mathbf u \in[0,1]^d} \bigl\llvert \hat
\Cb _n^{(m)}(0,1,\mathbf u)\bigr\rrvert, \qquad m \in\{1,\ldots,M
\}.
\end{equation}
The\vspace*{1pt} dimension $d$ was fixed to two, and the integrals and the suprema
in~(\ref{eqSn}) and~(\ref{eqTn}), respectively, were computed
approximately using a fine grid on $(0,1)^2$ of 400 uniformly spaced points.

Four data generating models were considered. 
The first one is a simple AR1 model. Let $\mathbf U_i$, $i \in\{-100,\ldots,0,\ldots,n\}$, be a bivariate i.i.d. sample from a copula $C$. Then,
set $\bolds\epsilon_i = (\Phi^{-1}(U_{i1}),\Phi^{-1}(U_{i2}))$, where
$\Phi$ is the c.d.f. of the standard normal distribution, and $\mathbf
X_{-100} = \bolds\epsilon_{-100}$. Finally, for any $j \in\{1,2\}$ and
$i \in\{-99,\ldots,0,\ldots,n\}$, compute recursively
{\renewcommand{\theequation}{AR1}
%
\begin{equation}\label{eq.AR1}
\label{eqar1} X_{ij} = 0.5
X_{i-1,j} + \epsilon_{ij}.
\end{equation}}%
The second and third data generating models are related to the
nonlinear autoregressive (NAR) model used in Paparoditis and Politis
(\cite{PapPol01}, Section~3.3), and to the exponential autoregressive
(EXPAR) model considered in Auestad and Tj{\o}stheim \cite{AueTjo90} and Paparoditis and Politis
(\cite{PapPol01}, Section~3.3). The sample $\mathbf X_1,\ldots,\mathbf X_n$ is
generated as previously with~(\ref{eqar1}) replaced by
{\renewcommand{\theequation}{NAR}
%
\begin{equation}\label{eq.NAR}
X_{ij} = 0.6 \sin(
X_{i-1,j} ) + \epsilon_{ij}
\end{equation}}%
and
{\renewcommand{\theequation}{EXPAR}
%
\begin{equation}
\label{eq.EXPAR}
X_{ij} = \bigl\{ 0.8 - 1.1 \exp
\bigl( - 50 X_{i-1,j}^2 \bigr) \bigr\} X_{i-1,j} + 0.1
\epsilon_{ij},
\end{equation}}%
respectively. The fourth and last data generating model is the
bivariate GARCH-like model considered in B{\"u}cher and Ruppert \cite
{BucRup13}. The sample of innovations is defined as for the models
above. In addition, for any $j \in\{1,2\}$, let $\sigma_{-100,j} =
\sqrt{\omega_j/(1-\alpha_j-\beta_j)}$ where $\omega_j$, $\alpha
_j$ and $\beta_j$ are usual $\operatorname{GARCH}(1,1)$ parameters whose values will be
set below, and, for any $j \in\{1,2\}$ and $i \in\{-99,\ldots,0,\ldots,n\}$, compute recursively
{\renewcommand{\theequation}{GARCH}
%
\begin{equation}
\label{eq.GARCH}
\sigma_{ij}^2 =
\omega_j + \beta_j \sigma_{i-1,j}^2
+ \alpha_j \epsilon_{i-1,j}^2 \quad\mbox{and}
\quad X_{ij} = \sigma_{ij} \epsilon_{ij}.
\end{equation}}\setcounter{equation}{2}%
Following B{\"u}cher and Ruppert \cite{BucRup13}, we take $(\omega
_1,\beta_1,\alpha_1)=(0.012,0.919,0.072)$ and $(\omega_2,\beta
_2,  \alpha_2)=(0.037,0.868,0.115)$. The latter values were estimated by
Jondeau, Poon and Rockinger \cite{JonPooRoc07} from SP500 and DAX daily logreturns, respectively.

The other factors of the experiments are as follows. Four different
copulas were considered: Clayton copulas with parameter values 1 and 4,
respectively, and Gumbel--Hougaard copulas with parameter value 1.5 and
3, respectively. The lower (resp., higher) parameter values correspond
to a Kendall's tau of $1/3$ (resp., $2/3$), that is, to mild (resp., strong) dependence. Notice that the Clayton copula is lower-tail
dependent while the Gumbel--Hougaard is upper-tail dependent
(see, e.g., McNeil, Frey and Embrechts \cite{McNFreEmb05}, Chapter~5).
The values 100, 200 and 400 were considered for $n$.

\begin{figure}

\includegraphics{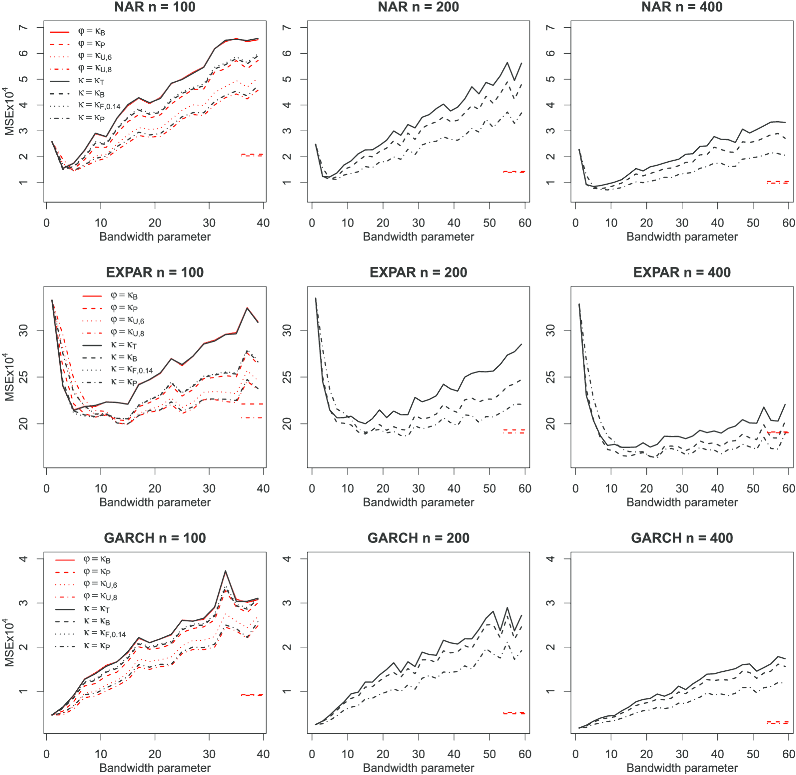}

\caption{For\vspace*{1pt} various choices of the function $\kappa/\varphi$ (see
Section~\protect\ref{secgenmult}), empirical MSE${}\times10^4$ of
the estimator $S_n^{(\lfloor0.95M \rfloor:  M)}$ with $M=2500$ versus
the bandwidth parameter $\ell_n$ under the NAR, EXPAR and GARCH data
generating scenarios with $C$ being the Gumbel--Hougaard copula with
parameter 1.5. The line segments in the lower-right corners of the
graphs correspond to the empirical MSEs of the estimator with estimated
bandwidth parameter following the procedure described in
Section~\protect\ref{secband}. The line styles of the segments
correspond to the choice of $\varphi$.}\label{figmse}
\end{figure}

We report the results of the experiments very partially (additional
results are available
in the supplementary material, see B\"ucher and Kojadinovic \cite
{BucKoj14}) and when based on the estimators of the partial derivatives
given\vspace*{1pt} in~(\ref{eqpdestRS}). Figure~\ref{figmse} displays the
empirical MSE of the estimator $S_n^{(\lfloor pM \rfloor:  M)}$ of the
quantile of order $p=0.95$ of $S_n$ versus the bandwidth parameter
$\ell_n$ for the different choices of $\kappa/\varphi$ mentioned
in Section~\ref{secgenmult}. The top (resp., middle, bottom) line of
graphs was obtained from datasets generated under the NAR (resp., EXPAR, GARCH) scenario with $C$ being the Gumbel--Hougaard copula with
parameter value~1.5. The line segments in the lower-right corners of
the graphs correspond\vspace*{1pt} to the empirical MSEs of the estimator
$S_n^{(\lfloor0.95M \rfloor:  M)}$ based on the estimated bandwidth
$\hat\ell_n^{\mathrm{opt}}$ computed as explained in Section~\ref{secband}.
The line styles of the segments correspond to the choice of $\varphi$.
The results for the AR1 scenario being very similar to those for the
NAR scenario are not reported. Similar looking graphs were obtained for
the other three copulas used in the simulations and when replacing the
Cram\'er--von Mises statistics by the Kolmogorov--Smirnov statistics
defined in~(\ref{eqTn}). In a related manner, the shapes of the
graphs were not too much affected by the value $p$ of the quantile
order: the empirical MSEs were smaller for $p < 0.95$ and higher for
$p=0.99$. Figures analogue to Figure~\ref{figmse} for other values of
$p$ and/or for the Kolmogorov--Smirnov statistic $T_n$ can actually be
found in Section~G of the supplementary material (B\"ucher
and Kojadinovic \cite{BucKoj14}).

The black (resp., red) curves in the first column of panels of
Figure~\ref{figmse} were obtained for dependent multiplier sequences
generated from initial standard normal i.i.d. sequences using the
moving average (resp., covariance matrix) approach described in
Section~\ref{subsecma} (resp., Section~\ref{subseccov}). The
functions $\kappa_T$, $\kappa_B$, $\kappa_{F,0.14}$ and $\kappa_P$
were considered for $\kappa$ in the case of the moving average
approach, while the function $\varphi$ in the covariance matrix
approach was successively taken equal to $\kappa_B$, $\kappa_P$,
$\kappa_{U,6}$ and $\kappa_{U,8}$. Looking at the graphs for $n=100$,
we see that, when the functions $\kappa$ and $\varphi$ are chosen to
match in the sense of Remark~\ref{remrel}, the resulting empirical
MSEs are very close. For that reason, to facilitate reading of the
plots, only the curves obtained with the moving average approach and
$\kappa\in\{\kappa_T, \kappa_B, \kappa_P\}$ are plotted when $n
\in\{200,400\}$. As it can be seen, for the NAR and EXPAR scenarios,
the empirical MSEs tend to decrease first with~$\ell_n$, reach a
minimum, and increase again. It is not the case for the GARCH setting
for which it seems that $\ell_n=1$ always leads to the smallest MSE.
In other words, the use of the dependent multiplier bootstrap does not
seem necessary in that context as the usual i.i.d. multiplier of R{\'e}millard and
Scaillet \cite{RemSca09} provides the best results. This might be due
to the fact that in this setting the contributions of the lagged
covariances to the long-run variance of the empirical process are very
small. Looking again at the graphs for the NAR and EXPAR settings, we
see that the smallest MSEs are reached by choosing $\kappa=\kappa
_P/\varphi=\kappa_{U,8}$, which is in accordance with
Proposition~\ref{propvar} which states that, asymptotically, kernels
with the smallest integral lead to the lowest variance. Another
observation is that, unlike what was expected by Shao (\cite{Sha10},
Remark 2.1) in the case of the mean as statistic of interest, the
choice $\kappa=\kappa_{F,0.14}$ did not lead to better results than
the choice $\kappa=\kappa_P$. Finally, let us comment on the
empirical MSEs of the estimator $S_n^{(\lfloor0.95M \rfloor:  M)}$
based on the estimated bandwidth $\hat\ell_n^{\mathrm{opt}}$ computed as
explained in Section~\ref{secband}. As it can be seen from the line
segments in the lower-right corners of the graphs, the achieved
empirical MSEs decrease with $n$ and are, overall, reasonably close to
the lowest observed MSE. Considering all the available results
(see Section~G of the supplementary material B\"ucher and Kojadinovic
\cite{BucKoj14}, for additional figures), the choice $\varphi= \kappa
_{U,8}$ appears to lead to a slightly lower MSE, overall, when $n=100$.
For $n \in\{200,400\}$, the choices $\varphi= \kappa_P$ and $\varphi
= \kappa_{U,8}$ do not seem to lead to differences of practical interest.

%

In view of the small differences between the moving average and
covariance matrix approaches for generating dependent multipliers
(black versus red curves in the first column of graphs of Figure~\ref
{figmse}), we suggest to use the former which is faster and more
stable numerically as it does not require the computation of the square
root of a large covariance matrix.

Before discussing further the estimation of $\ell_n$ using the results
of Section~\ref{secband}, let us mention an observation of practical
interest. Working with the same random seed, we replicated the
experiments described above using the two alternative definitions of
the partial derivative estimators mentioned below~(\ref{eqpdestRS}).
To our surprise, the best results, overall, were obtained with the
proposal of R{\'e}millard and
Scaillet \cite{RemSca09} given in~(\ref{eqpdestRS}), although the
differences seem too small to be of practical interest.

%
\begin{figure}

\includegraphics{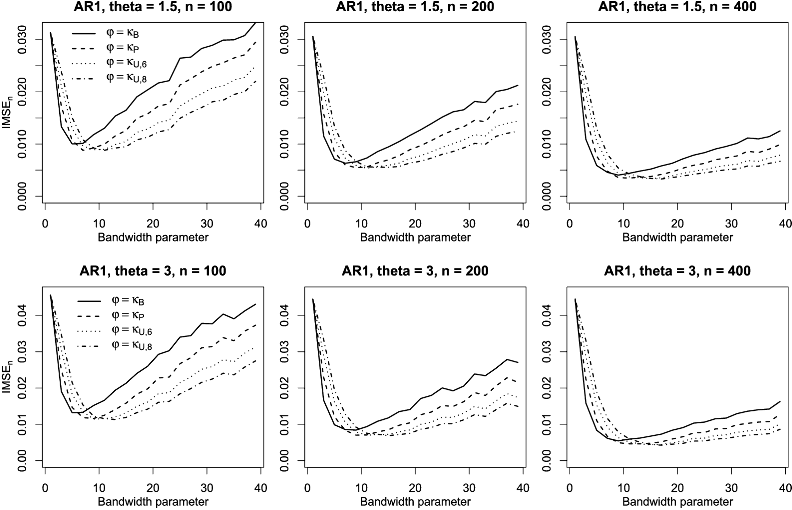}

\caption{For several choices of the function $\varphi$, $\IMSE_n$
defined in~\protect(\ref{eqimse}), computed approximately using a
grid of 25 uniformly spaced points on $(0,1)^2$ and $1000$ samples
versus the bandwidth parameter $\ell_n$ under the AR1 data generating
scenario with~$C$ being the Gumbel--Hougaard copula with parameter 1.5
(top row) and parameter~3 (bottom row).}\vspace*{10pt}
\label{figimse}
\end{figure}

We end this section with a more direct empirical investigation of the
estimator $\hat\ell_n^{\mathrm{opt}}$ of $\ell_n^{\mathrm{opt}}$ (see~(\ref
{eqlnopt}) and Section~\ref{secband}). We report an experiment based
on the AR1 model which will serve as a benchmark for judging about the
performance of $\hat\ell_n^{\mathrm{opt}}$. The setting is the following: a
grid $\{\mathbf u_i\}_{i \in\{1,\ldots,g\}}$ of $g=25$ points uniformly
spaced over $(0,1)^2$ was created, and $\sigma_C(\mathbf u_i, \mathbf u_j)$
was accurately estimated for all $i,j \in\{1,\ldots,g\}$ from $10^5$
samples of size $1000$ generated under the AR1 model described
previously. The latter estimation was carried out as follows: given a
sample $\mathbf X_1,\ldots,\mathbf X_n$ generated from the AR1 model, the
marginally standard uniform sample $\mathbf U_1,\ldots,\mathbf U_n$ was formed
using the fact that the marginal c.d.f.s of the $\mathbf X_i$ are centered
normal with variance $1/(1-0.5^2)$ in this case; this enabled us to
compute $\tilde\B_n(1,\cdot)$ at the grid points, where $\tilde\B
_n$ is defined in~(\ref{eqseqep}); for any $i,j \in\{1,\ldots,g\}$,
$\sigma_C(\mathbf u_i, \mathbf u_j)$ was finally accurately estimated as the
sample covariance of $10^5$ independent realizations of $(\tilde\B
_n(1,\mathbf u_i),\tilde\B_n(1,\mathbf u_j))$.

%
\begin{table}
\tabcolsep=0pt
\caption{Mean and standard deviation of 1000 estimates of $\ell
_n^{\mathrm{opt}}$, defined in~\protect(\ref{eqlnopt}), computed as explained
in Section~\protect\ref{secband} from 1000 samples generated from
the AR1 model in which $C$ is the Gumbel--Hougaard copula with
parameter $\theta$. The computations were carried out for the choices
$\varphi=\kappa_P$ and $\varphi=\kappa_{U,8}$}\label{tabbandBn}
\begin{tabular*}{\tablewidth}{@{\extracolsep{\fill}}@{}llllll@{}}
\hline
& & \multicolumn{2}{c}{$\varphi=\kappa_P$} & \multicolumn{2}{c}{$\varphi=\kappa_{U,8}$} \\ [-4pt]
& & \multicolumn{2}{c}{\hrulefill} & \multicolumn{2}{c}{\hrulefill}\\
$\theta$ & $n$ & Mean & Std. & Mean & Std. \\
\hline
1.5 & 100 & \phantom{0}8.93 & 3.85 & 12.41 & 5.92 \\
& 200 & 10.67 & 4.05 & 14.74 & 5.15 \\
& 400 & 12.81 & 3.94 & 17.73 & 4.99
\\[3pt]
3.0 & 100 & \phantom{0}9.11 & 5.18 & 12.75 & 8.13 \\
& 200 & 10.64 & 4.08 & 14.69 & 5.74 \\
& 400 & 12.77 & 3.94 & 17.66 & 5.31 \\
\hline
\end{tabular*}
\end{table}

Next, for $n \in\{100,200,400\}$ and $\ell_n \in\{1,3,\ldots,39\}$,
$\IMSE_n$ defined in~(\ref{eqimse}) was approximated as follows:
$1000$ samples $\mathbf X_1,\ldots,\mathbf X_n$ were generated under the AR1
model, and, for each sample, the processes $\hat\B_n^{(1)}(1,\cdot
),\ldots,\hat\B_n^{(M)}(1,\cdot)$ with $M=1000$ were evaluated at
the grid points, with $\hat\B_n^{(m)}$ defined in~(\ref{eqhatBnm});
computing sample covariances, this allowed us to obtain $1000$
bootstrap estimates of $\sigma_C(\mathbf u_i, \mathbf u_j)$ for all $i,j \in\{
1,\ldots,g\}$, from which we approximated $\IMSE_n$. The results are
represented in the graphs of Figure~\ref{figimse} for the previously
considered choices of the function $\varphi$. The top (resp., bottom)
row of graphs was obtained when $C$ in the AR1 data generating scenario
is the Gumbel--Hougaard copula with parameter 1.5 (resp., 3).

The procedure described in Section~\ref{secband} was finally used to
obtain 1000 estimates of $\ell_n^{\mathrm{opt}}$ under the AR1 model based on
the Gumbel--Hougaard copula with parameter $\theta$, for $n \in\{
100,200,400\}$, $\varphi\in\{\kappa_P,\kappa_{U,8}\}$ and $\theta
\in\{1.5,3\}$. The mean and standard deviation of the estimates are
reported in Table~\ref{tabbandBn}. A comparison with Figure~\ref
{figimse} reveals that the procedure described in Section~\ref
{secband} for estimating $\ell_n^{\mathrm{opt}}$ gives surprisingly good
results on average for the experiment at hand. Another observation is
that the estimates do not seem much affected by the value of $\theta$,
that is, the strength of the dependence.


\begin{appendix}
\section{Proof of Theorem~\texorpdfstring{\protect\ref{teomultCLT}}{2.1}}\label{proofuncondBuh}

The proof of Theorem~\ref{teomultCLT} is based on three lemmas. The
first lemma establishes weak convergence of the finite-dimensional
distributions, while the second and third lemmas concern asymptotic tightness.


The following result can be proved using a well-known blocking
technique (see, e.g., Dehling and Philipp \cite{DehPhi02}, page 31).
Its proof is given in the supplementary material
(B\"ucher and Kojadinovic \cite{BucKoj14}).

%
\begin{lem}[(Finite-dimensional convergence)]
\label{lemfidi}
Assume that $\ell_n = \mathrm{O}(n^{1/2 - \varepsilon})$ for some $0 <
\varepsilon< 1/2$ and that $(\mathbf U_i)_{i \in\Z}$ is a strictly
stationary sequence whose strong mixing coefficients satisfy $\alpha_r
= \mathrm{O}(r^{-a})$, $a > 2$. Then, the finite-dimensional distributions of
$(\tilde\B_n, \tilde\B_n^{(1)},\ldots, \tilde\B_n^{(M)} )$
converge weakly to those of $(\B_C, \B_C^{(1)},\ldots, \B_C^{(M)} )$.
\end{lem}


Regarding the tightness, let us first extend $\tilde\B_n^{(m)}$, $m
\in\{1,\ldots,M\}$, to blocks in $[0,1]^{d+1}$ in the spirit of Bickel
and Wichura \cite{BicWic71}. For any $(s,t] \subset[0,1]$ and $A =
(u_1,v_1] \times\cdots\times(u_d,v_d] \subset[0,1]^d$, we define
$\tilde\B_n^{(m)}((s,t] \times A)$ to be
\[
\tilde\B_n^{(m)}\bigl((s,t] \times A\bigr) = \frac{1}{\sqrt{n}} \sum
_{i=\lfloor ns \rfloor + 1}^{\lfloor nt \rfloor} \xi_{i,n}^{(m)}
\bigl[ \1 (\mathbf U_{i} \in A) - \nu(A) \bigr], %
\]
where
\begin{eqnarray*}
\nu(A) &=& \Pr(\mathbf U_1 \in A)
\\
&=& \sum_{(\epsilon_1,\ldots,\epsilon_d)
\in\{0,1\}^d}
(-1)^{\sum_{i=1}^d \epsilon_i} C\bigl\{(1-\epsilon_1) v_1 +
\epsilon_1 u_1,\ldots,(1-\epsilon_d)
v_d + \epsilon_d u_d \bigr\}. %
\end{eqnarray*}

In the next two lemmas, the sequences $(\xi_{i,n}^{(m)})_{i \in\Z}$
are only assumed to satisfy \textup{(M1)} with $\Ex[ \{ \xi
_{0,n}^{(m)} \}^2 ] > 0$ not necessarily equal to one.

%
\begin{lem}[(Moment inequality)]\label{lem4mom}
Assume that $(\mathbf U_i)_{i \in\Z}$ is a strictly stationary sequence
whose strong mixing coefficients satisfy $\alpha_r = \mathrm{O}(r^{-a})$, $a >
6$. Then,\vspace*{1pt} for any $m \in\{1,\ldots,M\}$, $q \in(2a/(a-3), 4)$, $(s,t]
\subset[0,1]$ and $A = (u_1,v_1] \times\cdots\times(u_d,v_d]
\subset[0,1]^d$, we have
\[
\Ex \bigl[ \bigl\{ \tilde\B_n^{(m)}\bigl((s,t] \times
A\bigr) \bigr\}^4 \bigr] \leq \kappa\bigl[ \lambda_n(s,t)^2
\bigl\{\nu(A) \bigr\}^{4/q} + n^{-1} \lambda_n(s,t)
\bigl\{\nu(A)\bigr\}^{2/q} \bigr], %
\]
where $\kappa>0$ is a constant.
\end{lem}

\begin{pf}
The proof is similar to that of Lemma 3.22 in Dehling and Philipp \cite
{DehPhi02}. Fix $m \in\{1,\ldots,M\}$. For any $i \in\Z$, let $Y_i =
\1(\mathbf U_i \in A) - \nu(A)$. Then,
%
\begin{eqnarray}
\label{eqsum1}
&&\Ex \bigl[ \bigl\{ \tilde\B_n^{(m)}\bigl(
(s,t] \times A\bigr) \bigr\}^4 \bigr]\nonumber
\\
&&\quad  = \frac{1}{n^2} \sum
_{i_1, i_2, i_3, i_4=\lfloor ns \rfloor
+1}^{\lfloor nt \rfloor} \Ex\bigl[\xi_{i_1,n}^{(m)}
\xi_{i_2,n}^{(m)} \xi _{i_3,n}^{(m)}
\xi_{i_4,n}^{(m)}\bigr] \Ex[Y_{i_1}Y_{i_2}
Y_{i_3}Y_{i_4}]
\\
&&\quad \le\frac{4!\lambda_n(s,t)}{ n} \mathop{\sum_{0\le i,j,k \le\lfloor nt
\rfloor-\lfloor ns \rfloor-1}}_{i+j+k \le\lfloor nt \rfloor
-\lfloor ns \rfloor-1} \bigl\llvert \Ex
\bigl[ \xi_{0,n}^{(m)}\xi_{i,n}^{(m)} \xi
_{i+j,n}^{(m)} \xi_{i+j+k,n}^{(m)}\bigr]
\Ex[Y_{0}Y_{i} Y_{i+j}Y_{i+j+k}] \bigr
\rrvert.\nonumber
\end{eqnarray}
On one hand, $\llvert  \Ex[\xi_{0,n}^{(m)}\xi_{i,n}^{(m)} \xi_{i+j,n}^{(m)}
\xi_{i+j+k,n}^{(m)}] \rrvert  \le\Ex[\{ \xi_{0,n}^{(m)}\}^4]$. On the other
hand, by Lemma 3.11 of Dehling and Philipp \cite{DehPhi02}, for any $q
\in(2a/(a-3), 4)$ and $p \in(2,a/3)$ such that $1/p+2/q=1$, we have
\begin{eqnarray*}
\Ex\bigl[Y_{0} (Y_{i} Y_{i+j}Y_{i+j+k})
\bigr] &\le&10 \alpha_i^{1/p} \llVert Y_0
\rrVert _q \llVert Y_iY_{i+j}Y_{i+j+k}
\rrVert _q \le10 \alpha_i^{1/p} \llVert
Y_0\rrVert _q^2,
\\
\Ex\bigl[(Y_{0} Y_{i} Y_{i+j})Y_{i+j+k}
\bigr] &\le&10 \alpha_k^{1/p} \llVert Y_0
\rrVert _q^2
\end{eqnarray*}
and
\begin{eqnarray*}
\bigl\llvert \Ex\bigl[(Y_{0} Y_{i}) (
Y_{i+j}Y_{i+j+k}) \bigr] \bigr\rrvert &\le&\bigl\llvert
\Ex[Y_{0} Y_{i} ] \Ex [ Y_{i+j}Y_{i+j+k}
] \bigr\rrvert + 10 \alpha_j^{1/p} \llVert
Y_0Y_i\rrVert _q \llVert
Y_{i+j}Y_{i+j+k} \rrVert _q
\\
&\le& 100 \alpha_i^{1/p} \alpha_k^{1/p}
\llVert Y_0\rrVert _q^4 + 10 \alpha
_j^{1/p} \llVert Y_0\rrVert
_q^2.
\end{eqnarray*}
Proceeding as in Lemma 3.22 of Dehling and Philipp \cite{DehPhi02}, we
split the sum on the right of~(\ref{eqsum1}) into three sums
according to which of the indices $i,j,k$ is the largest. Combining
this decomposition with the three previous inequalities, we obtain
\begin{eqnarray*}
&& \Ex \bigl[ \bigl\llvert \tilde\B_n^{(m)}\bigl((s,t]
\times A\bigr) \bigr\rrvert ^4\bigr]
\\
&&\quad \le\frac{24 \Ex[\{\xi_{0,n}^{(m)}\}^4] \lambda_n(s,t)}{n}
\\
&&\qquad{}\times \Biggl\{ 100
\llVert Y_0\rrVert _q^4 \sum
_{j=0}^{\lfloor nt \rfloor-\lfloor ns \rfloor
-1} \sum_{i,k \le j}
\alpha_i^{1/p} \alpha_k^{1/p} + 30
\llVert Y_0\rrVert _q^2 \sum
_{i=0}^{\lfloor nt \rfloor-\lfloor ns \rfloor
-1} \sum_{j,k\le i}
\alpha_i^{1/p} \Biggr\}.
\end{eqnarray*}
Observing that $\sum_{i=1}^\infty\alpha_i^{1/p} < \infty$ and $\sum_{i=1}^\infty i^2 \alpha_i^{1/p}< \infty$ (note that $p<a/3$ by
construction), we can bound the expression on the right of the previous
inequality by
\[
\kappa \bigl\{ \lambda_{n}(s,t)^2\llVert
Y_0\rrVert _q^4 + n^{-1} \lambda
_n(s,t) \llVert Y_0\rrVert _q^2
\bigr\},
\]
where $\kappa> 0$ is a constant depending on the mixing coefficients
and $\Ex[\{\xi_{0,n}^{(m)}\}^4]$. Finally, since \mbox{$q>2$} by
construction, the assertion follows from the fact that $\Ex[\llvert  Y_0 \rrvert  ^q]
\le\Ex[Y_0^2] =\nu(A)-\nu(A)^2 \le\nu(A)$.
\end{pf}


Let us introduce additional notation. For any $\delta\geq0$, $T
\subset[0,1]^{d+1}$ and $f \in\ell^\infty([0,1]^{d+1})$, let
\[
w_\delta(f,T) = \mathop{\sup_{x, y \in T}}_{\llVert  x-y\rrVert  _1 \le\delta} \bigl\llvert f(x) - f(y)
\bigr\rrvert, %
\]
where $\llVert  \cdot \rrVert  _1$ denotes the 1-norm.

%
\begin{lem}[(Asymptotic equicontinuity)]
\label{lemasymequi}
Assume that $(\mathbf U_i)_{i \in\Z}$ is a strictly stationary sequence
whose strong mixing coefficients satisfy $\alpha_r = \mathrm{O}(r^{-a})$, $a >
3 + 3d/2$. Then, for any $m \in\{1,\ldots,M\}$, $\tilde\B_n^{(m)}$
is asymptotically uniformly $\llVert  \cdot \rrVert  _1$-equicontinuous in
probability, that is, for any $\varepsilon> 0$,
\[
\lim_{\delta\downarrow0} \limsup_{n\to\infty} \Pr\bigl\{
w_{\delta} \bigl(\tilde\B_n^{(m)},
[0,1]^{d+1} \bigr) > \varepsilon\bigr\} = 0. %
\]
\end{lem}

\begin{pf}
Fix $m \in\{1,\ldots,M\}$. Let $K > 0$ be a constant and let us first
assume that, for any $n \geq1$ and $i \in\{1,\ldots,n\}$, $\xi
_{i,n}^{(m)} \geq-K$. Then, let $Z_{i,n}^{(m)} = \xi_{i,n}^{(m)} + K
\geq0$. Furthermore, let $\gamma\in(0,1/2]$ be a real parameter to
be chosen later, and define
\[
I_n=\{i/n:   i=0,\ldots,n\}, \qquad I_{n,\gamma}=\bigl\{i/\bigl
\lfloor n^{1/2+\gamma} \bigr\rfloor:   i=0,\ldots, \bigl\lfloor n^{1/2+\gamma}
\bigr\rfloor \bigr\},
\]
and $T_n = I_n \times I_{n,\gamma}^d$. Also, for any $s \in[0,1]$,
let $\underline{s}=\lfloor sn \rfloor/n$ and $\bar{s} = \lceil sn
\rceil/n$; clearly, $\underline{s},\bar{s} \in I_n$ and are such
that $\underline{s} \leq s \leq\bar{s}$ and $\bar{s} - \underline
{s} \leq1/n$. Similarly, for any $u \in[0,1]$, let $\underline
{u}_\gamma, \bar{u}_\gamma\in I_{n,\gamma}$ such that $\underline
{u}_\gamma\leq u \leq\bar{u}_\gamma$ and $\bar{u}_\gamma-
\underline{u}_\gamma\leq1/\lfloor n^{1/2+\gamma} \rfloor$. Then,
for any $\mathbf u \in[0,1]^d$, we define $\underline{\mathbf u}_\gamma\in
I_{n,\gamma}^d$ (resp., $\bar{\mathbf u}_\gamma\in I_{n,\gamma}^d$) as
$\underline{\mathbf u}_\gamma= (\underline{u}_{1,\gamma},\ldots,\underline{u}_{d,\gamma})$ (resp., $\bar{\mathbf u}_\gamma= (\bar
{u}_{1,\gamma},\ldots,\bar{u}_{d,\gamma})$).

Now, for any $(s,\mathbf u) \in[0,1]^{d+1}$,
\begin{eqnarray*}
\tilde\B_n^{(m)}(s,\mathbf u) - \tilde\B_n^{(m)}(s,
\underline{\mathbf u}_\gamma)
&\le&\frac{1}{\sqrt{n}} \sum
_{i=1}^{\lfloor ns \rfloor} Z_{i,n}^{(m)} \bigl\{
\1(\mathbf U_i \le\bar{\mathbf u}_\gamma) - \1(\mathbf U_i
\le \underline{\mathbf u}_\gamma) \bigr\}
\\
&&{} + \sqrt n K \bigl\{ C(\bar{\mathbf u}_\gamma ) - C(\underline{\mathbf
u}_\gamma) \bigr\}.
\end{eqnarray*}
Thus,
\begin{eqnarray*}
\tilde\B_n^{(m)}(s,\mathbf u) - \tilde\B_n^{(m)}(s,
\underline{\mathbf u}_\gamma)
&\le&\tilde\B_n^{(m)}(s,\bar{
\mathbf u}_\gamma) -\tilde\B _n^{(m)}(s,\underline{\mathbf
u}_\gamma)
+ K \bigl\{ \tilde\B_n(s,\bar {\mathbf u}_\gamma) - \tilde
\B_n(s,\underline{\mathbf u}_\gamma) \bigr\}
\\
&&{} + \Biggl( \sqrt n K +
\frac{1}{\sqrt{n}} \sum_{i=1}^{\lfloor ns
\rfloor}
Z_{i,n}^{(m)} \Biggr) \bigl\{ C(\bar{\mathbf u}_\gamma) -
C(\underline{\mathbf u}_\gamma) \bigr\},
\end{eqnarray*}
and therefore
\begin{eqnarray*}
\tilde\B_n^{(m)}(s,\mathbf u) - \tilde\B_n^{(m)}(s,
\underline{\mathbf u}_\gamma) &\le& \bigl\llvert \tilde\B_n^{(m)}(s,
\bar{\mathbf u}_\gamma) -\tilde\B_n^{(m)}(s,\underline{
\mathbf u}_\gamma) \bigr\rrvert
+ K \bigl\llvert \tilde\B_n(s,\bar{\mathbf u}_\gamma) - \tilde
\B_n(s,\underline{\mathbf u}_\gamma) \bigr\rrvert
\\
&&{}+ d
\bigl(n^{\gamma} - 1\bigr)^{-1} \Bigl( K + \max
_{1 \leq i
\leq n} \bigl\llvert Z_{i,n}^{(m)}\bigr
\rrvert \Bigr),
\end{eqnarray*}
using the fact that $C$ satisfies the Lipschitz condition
%
\begin{equation}
\label{eqlipcond} \bigl\llvert C(\mathbf u) - C(\mathbf v)\bigr\rrvert \leq\llVert \mathbf u -
\mathbf v \rrVert _1\qquad\forall \mathbf u, \mathbf v \in[0,1]^d,
\end{equation}
and that $n^{1/2} (\lfloor n^{1/2 + \gamma} \rfloor)^{-1} \leq
(n^\gamma- 1)^{-1}$ for all $n \geq1$. Similarly, for any $(s,\mathbf u)
\in[0,1]^{d+1}$,
\begin{eqnarray*}
&& \tilde\B_n^{(m)}(s,\underline{\mathbf u}_\gamma) -
\tilde\B _n^{(m)}(s, \mathbf u)
\\
&&\quad \le\frac{1}{\sqrt{n}} \sum
_{i=1}^{\lfloor ns
\rfloor} Z_{i,n}^{(m)}
\bigl\{ C(\bar{\mathbf u}_\gamma)-C(\underline{\mathbf u}_\gamma) \bigr\}
+ \frac{K}{\sqrt{n}} \sum_{i=1}^{\lfloor ns
\rfloor} \bigl
\{ \1(\mathbf U_i \le\bar{\mathbf u}_\gamma) - \1(\mathbf
U_i \le \underline{\mathbf u}_\gamma) \bigr\}
\\
&&\quad \le d\bigl(n^{\gamma} - 1\bigr)^{-1} \Bigl(K + \max
_{1 \leq i \leq n} \bigl\llvert Z_{i,n}^{(m)}\bigr
\rrvert \Bigr) + K \bigl\llvert \tilde\B_n(s,\bar{\mathbf
u}_\gamma) - \tilde\B_n(s,\underline{\mathbf u}_\gamma)
\bigr\rrvert.
\end{eqnarray*}
Hence, for any $(s,\mathbf u) \in[0,1]^{d+1}$, we have that
%
\begin{eqnarray}
\label{eqchaining}
&& \bigl\llvert \tilde\B_n^{(m)}(s, \mathbf u) -
\tilde\B_n^{(m)}(s,\underline {\mathbf u}_\gamma)\bigr
\rrvert\nonumber
\\
&&\quad  \le \bigl\llvert \tilde\B_n^{(m)}(s,\bar{\mathbf
u}_\gamma) -\tilde\B_n^{(m)}(s,\underline{\mathbf
u}_\gamma) \bigr\rrvert+ K \bigl\llvert \tilde\B_n(s,\bar{\mathbf u}_\gamma) - \tilde\B
_n(s,\underline{\mathbf u}_\gamma) \bigr\rrvert
\\
&&\qquad  + d
\bigl(n^{\gamma} - 1\bigr)^{-1} \Bigl( K + \max
_{1 \leq i \leq n}\bigl\llvert Z_{i,n}^{(m)}\bigr\rrvert
\Bigr).\nonumber
\end{eqnarray}
Then, noticing that, for any $s \in[0,1]$, $\tilde\B_n^{(m)}(s,\cdot
) = \tilde\B_n^{(m)}(\underline{s},\cdot)$, and applying~(\ref
{eqchaining}) to the first and the third summand on the right-hand
side of the decomposition
\begin{eqnarray*}
\tilde\B_n^{(m)}(s,\mathbf u) -\tilde\B_n^{(m)}(t,
\mathbf v)
&=& \bigl\{ \tilde \B_n^{(m)}(\underline{s}, \mathbf u) -
\tilde\B_n^{(m)}(\underline{s},\underline{\mathbf
u}_\gamma) \bigr\}
+ \bigl\{ \tilde\B_n^{(m)}(\underline{s},\underline{\mathbf
u}_\gamma) - \tilde\B_n^{(m)}(\underline{t},
\underline{\mathbf v}_\gamma) \bigr\}
\\
&&{}+ \bigl\{ \tilde\B_n^{(m)}(
\underline{t},\underline{\mathbf v}_\gamma) - \tilde \B_n^{(m)}(
\underline{t}, \mathbf v) \bigr\},
\end{eqnarray*}
we obtain that, for any $\delta> 0$,
\begin{eqnarray*}
w_\delta\bigl(\tilde\B_n^{(m)},
[0,1]^{d+1} \bigr) &\le& 3 w_{\delta+(d+1)/
\lfloor n^{1/2+\gamma} \rfloor} \bigl(\tilde
\B_n^{(m)}, T_n \bigr) + 2 K w_{\delta+d/\lfloor n^{1/2+\gamma} \rfloor}
\bigl(\tilde\B_n, [0,1]^{d+1}\bigr)
\\
&&{}+ 2d\bigl(n^{\gamma} - 1\bigr)^{-1} \Bigl( K + \max
_{1 \leq i \leq n} \bigl\llvert Z_{i,n}^{(m)}\bigr
\rrvert \Bigr)
\\
&\le& 3 w_{2\delta} \bigl(\tilde\B_n^{(m)},
T_n \bigr) + 2 K w_{2\delta} \bigl(\tilde\B_n,
[0,1]^{d+1}\bigr)
\\
&&{}+ 2d\bigl(n^{\gamma} - 1\bigr)^{-1}
\Bigl( K + \max_{1
\leq i \leq n} \bigl\llvert Z_{i,n}^{(m)}
\bigr\rrvert \Bigr),
\end{eqnarray*}
for sufficiently large $n$. Now, from the previous inequality, for any
$\varepsilon> 0$,
\begin{eqnarray*}
\Pr\bigl\{w_\delta \bigl(\tilde\B_n^{(m)},
[0,1]^{d+1} \bigr)> \varepsilon\bigr\} &\le& \Pr\bigl\{ 3 w_{2\delta}
\bigl(\tilde\B_n^{(m)}, T_n \bigr)>\varepsilon/3
\bigr\}
\\
&&{}+ \Pr\bigl\{ 2 K w_{2\delta} \bigl(\tilde\B_n,
[0,1]^{d+1}\bigr) >\varepsilon/3\bigr\}
\\
&&{}+ \Pr\Bigl\{ 2d
\bigl(n^{\gamma} - 1\bigr)^{-1} \Bigl( K + \max
_{1 \leq i \leq n} \bigl\llvert Z_{i,n}^{(m)}\bigr
\rrvert \Bigr) > \varepsilon/3\Bigr\}.
\end{eqnarray*}
Since $a > 1$, we have from B\"ucher (\cite{Buc14}, Lemma~2) that
$\tilde\B_n$ is asymptotically uniformly \mbox{$\llVert  \cdot \rrVert
_1$-}equicontinuous in probability. This implies that the second term on
the right of the previous display converges to $0$ as $n \to\infty$
followed by $\delta\downarrow0$. The third term converges to zero
because $n^{-\gamma} \max_{1 \leq i \leq n} \llvert  Z_{i,n}^{(m)}\rrvert   \stackrel
{\Pr}{\to}0$. Indeed, for any $\eta> 0$ and $\nu> 1/\gamma\geq2$,
by Markov's inequality and~\textup{(M1)},
\[
\Pr\Bigl(n^{-\gamma} \max_{1 \leq i \leq n} \bigl\llvert
Z_{i,n}^{(m)}\bigr\rrvert > \eta\Bigr) \leq n \Pr\bigl( \bigl
\llvert Z_{1,n}^{(m)}\bigr\rrvert \geq\eta n^\gamma
\bigr) \leq\eta^{-\nu} n^{1 -
\gamma\nu} \sup_{n \geq1} \Ex
\bigl(\bigl\llvert Z_{1,n}^{(m)}\bigr\rrvert ^\nu
\bigr) \to0. %
\]
Thus, it remains to show that, for any $\varepsilon> 0$, $\lim_{\delta\downarrow0} \limsup_{n\to\infty} \Pr\{ w_{\delta}
(\tilde\B_n^{(m)}, T_n ) > \varepsilon\} = 0$, or equivalently
(see, e.g., van der Vaart and
Wellner \cite{vanWel96}, Problem 2.1.5) that, for any positive
sequence $\delta_n \downarrow0$, $\lim_{n\to\infty} \Pr\{
w_{\delta_n} (\tilde\B_n^{(m)}, T_n ) > \varepsilon\} = 0$. To do
so, we shall use Lemma~\ref{lem4mom} together with Lemma~2 of
Balacheff and Dupont
\cite{BalDup80} (see also Bickel and Wichura \cite{BicWic71}, Theorem~3 and the remarks on page 1665).

Recall that $\nu$ is the measure on $[0,1]^d$ corresponding to the
c.d.f. $C$, and let $\mu$ be a measure on $[0,1]^{d+1}$ defined by
$\mu=2\lambda\otimes\nu$, where $\lambda$ denotes the
one-dimensional Lebesgue measure. Next, for some real $q \in
(2a/(a-3),6a/(2a-3)) \subset(2,4)$, let $\beta= 2 - 2/q -3/a\in
(1,4/q)$. Furthermore, consider a non-empty set $(s,t] \times A= (s,t]
\times(u_1,v_1] \times\cdots\times(u_d,v_d]$ of $[0,1]^{d+1}$ whose
boundary points are all distinct and lie in $T_n$. Then, starting from
Lemma~\ref{lem4mom}, for any $q \in(2a/(a-3),6a/(2a-3)) \subset(2,4)$,
\begin{eqnarray*}
&& \Ex \bigl[ \bigl\{ \tilde\B_n^{(m)}\bigl((s,t] \times
A\bigr) \bigr\}^4 \bigr]
\\[-1pt]
&&\quad \le \kappa\bigl[ \lambda_n(s,t)^2
\bigl\{ \nu(A) \bigr\}^{4/q} + n^{-1} \lambda
_n(s,t)\bigl\{ \nu(A)\bigr\}^{2/q} \bigr]
\\[-1pt]
&&\quad \le\kappa\bigl[ \bigl\{ \lambda_n(s,t) \nu(A) \bigr
\}^{4/q} + n^{-1} \bigl\{ \lambda _n(s,t) \nu(A)
\bigr\}^{2/q}\bigr]
\\[-1pt]
&&\quad \le \kappa\mu\bigl((s,t] \times A\bigr)^\beta \bigl\{ \mu\bigl((s,t] \times
A\bigr)^{4/q - \beta} + n^{-1} \mu\bigl( (s,t] \times A\bigr)^{2/q-\beta} \bigr\}
\\[-1pt]
&&\quad \le \kappa\mu\bigl((s,t] \times A\bigr)^\beta \bigl\{ 2^{4/q - \beta} +
n^{-1} n^{-(1+d/2 +d\gamma)(2/q-\beta)} \bigr\}
\\[-1pt]
&&\quad = \kappa\mu\bigl((s,t] \times A\bigr)^\beta\bigl\{ 2^{4/q - \beta} +
n^{(\beta
-2/q)(1+d/2+d\gamma)-1} \bigr\}.
\end{eqnarray*}
Note that $\inf_{q > 2a/(a-3)} (\beta- 2/q) = 3/a$. Hence, because
$3/a < 2/(2+d)$ from the assumption on the mixing rate, it is possible
to choose $q \in(2a/(a-3),6a/(2a-3))$ and $\gamma> 0$ (the parameter
involved in the grid $I_{n,\gamma}^d$) small enough such that $\beta-
2/q < 2/(2+d+2d\gamma)$. For the aforementioned parameter choices,
$(\beta-2/q)(1+d/2+d\gamma)-1 < 0$, which implies that $n^{(\beta
-2/q)(1+d/2+d\gamma)-1} \leq1$ for all $n \geq1$.

With some abuse of notation consisting of incorporating the constant
$\{ \kappa( 2^{4/q - \beta} + 1) \}^{1/\beta}$ into the measure, we obtain
\[
\Ex \bigl[ \bigl\{ \tilde\B_n^{(m)}\bigl((s,t] \times
A\bigr) \bigr\}^4 \bigr] \le\mu \bigl((s,t] \times A\bigr)^\beta,
\]
which, by Markov's inequality, implies that, for any $\varepsilon> 0$,
\[
\Pr \bigl\{ \bigl\llvert \tilde\B_n^{(m)}\bigl((s,t] \times A\bigr)
\bigr\rrvert \ge\varepsilon \bigr\} \le\varepsilon^{-4} \mu\bigl((s,t] \times
A\bigr)^\beta.
\]
Now, let $\tilde\mu_n$ denote a finite measure on $T_n$ defined from
its values on the singletons $\{(s,\mathbf u)\}$ of $T_n$ as
\[
\tilde\mu_n\bigl(\bigl\{(s,\mathbf u)\bigr\}\bigr) = \cases{0, &\quad if
$s \wedge u_1 \wedge\cdots\wedge u_d = 0$,
\cr
\mu\bigl(
\bigl(s', s\bigr] \times\bigl(u_1', u_1\bigr]
\times\cdots\times\bigl(u_d', u_d\bigr]\bigr), &
\quad otherwise,}
\]
where $s'=\max\{ t \in I_n:   t < s\}$ and $u_j'=\max\{ u \in
I_{n,\gamma}:   u < u_j\}$ for all $j \in\{1,\ldots,d\}$. By
additivity of $\tilde\mu_n$, the previous estimation reads
\[
\Pr \bigl\{ \bigl\llvert \tilde\B_n^{(m)}\bigl((s,t] \times A\bigr)
\bigr\rrvert \ge\varepsilon \bigr\} \le\varepsilon^{-4} \tilde
\mu_n\bigl[ \bigl\{(s,t] \times A \bigr\} \cap T_n \bigr]
^\beta.
\]

We shall now conclude by an application of Lemma~2 of Balacheff and
Dupont \cite{BalDup80}. Consider a positive sequence $\delta_n
\downarrow0$, and let $\delta_n' \downarrow0$ such that, for any $n
\in\N$, $\delta_n' \in\{1/i:   i \in\N\}$ and $\delta_n' \ge\max
\{\delta_n, 1/\lfloor n^{1/2+\gamma} \rfloor\}$. Applying Lemma 2 of
Balacheff and Dupont \cite{BalDup80} (note that $1/\lfloor
n^{1/2+\gamma} \rfloor = \max\{1/n, 1/\lfloor n^{1/2+\gamma}
\rfloor\}$ is denoted by $\tau$ in the lemma) and using the fact that
$\llVert  \cdot \rrVert  _2 \leq\llVert  \cdot \rrVert  _1$, we obtain that, for any $\varepsilon
> 0$, there exists a constant $\lambda> 0$ depending on $\varepsilon
$, $\beta$ and $d$, such that
\begin{eqnarray*}
&& \Pr\bigl\{ w_{\delta_n} \bigl(\tilde\B_n^{(m)},
T_n \bigr) > \varepsilon\bigr\}
\\
&&\quad  \le \Pr\bigl\{ w_{\delta_n'}
\bigl(\tilde\B_n^{(m)}, T_n \bigr) > \varepsilon
\bigr\}
\\
&&\quad \le \lambda\tilde\mu_n(T_n)
\\
&&\qquad{}\times
\Bigl[ \max\Bigl\{ \mathop{\sup_{s,t \in I_n }}_{\llvert  s - t\rrvert   \leq3 \delta
_n'} \bigl\llvert \tilde
\mu_n\bigl( \{0,\ldots,s\} \times I_{n,\gamma}^d
\bigr) - \tilde \mu_n\bigl(\{0,\ldots,t\} \times I_{n,\gamma}^d
\bigr) \bigr\rrvert,
\\
&&\hspace*{58pt}\mathop{\sup_{u,v \in I_{n,\gamma}}}_{\llvert  u - v\rrvert   \leq3 \delta_n'} \bigl\llvert \tilde \mu_n
\bigl(I_n \times\{0,\ldots,u\} \times I_{n,\gamma}^{d-1}
\bigr) - \tilde \mu_n\bigl(I_n \times\{0,\ldots,v\} \times
I_{n,\gamma}^{d-1}\bigr) \bigr\rrvert,
\\
&&\hspace*{68pt} \dots,
\\
&&\hspace*{58pt}\mathop{\sup_{u,v \in I_{n,\gamma}}}_{\llvert  u - v\rrvert   \leq3 \delta_n'} \bigl\llvert \tilde \mu_n
\bigl(I_n \times I_{n,\gamma}^{d-1} \times\{0,\ldots,u\}
\bigr) - \tilde \mu_n\bigl(I_n \times I_{n,\gamma}^{d-1}
\times\{0,\ldots,v\}\bigr) \bigr\rrvert \Bigr\} \Bigr]^{\beta-1},
\end{eqnarray*}
which implies that,
\begin{eqnarray*}
&& \Pr\bigl\{ w_{\delta_n} \bigl(\tilde\B_n^{(m)},
T_n \bigr) > \varepsilon\bigr\}
\\
&&\quad  \le \lambda \mu\bigl([0,1]^{d+1}
\bigr)
\\
&&\qquad{} \times \Bigl[ \max\Bigl\{ \mathop{\sup_{s,t \in
[0,1] }}_{\llvert  s - t\rrvert   \leq3 \delta_n'} \bigl\llvert \mu
\bigl([0,s] \times[0,1]^d\bigr) - \mu\bigl([0,t] \times[0,1]^d
\bigr) \bigr\rrvert,
\\
&&\hspace*{58pt}\mathop{\sup_{u,v \in[0,1] }}_{\llvert  u - v\rrvert   \leq3 \delta_n'} \bigl\llvert \mu\bigl([0,1] \times[0,u]
\times[0,1]^{d-1}\bigr) - \mu\bigl([0,1] \times[0,v] \times
[0,1]^{d-1}\bigr) \bigr\rrvert,
\\
&&\hspace*{68pt} \dots,
\\
&&\hspace*{147pt} \mathop{\sup_{u,v \in[0,1] }}_{\llvert  u - v\rrvert   \leq3 \delta_n'} \bigl\llvert \mu\bigl([0,1]^d
\times[0,u]\bigr) - \mu\bigl([0,1]^d \times[0,v]\bigr) \bigr\rrvert
\Bigr\} \Bigr]^{\beta-1},
\end{eqnarray*}
which converges to $0$ by uniform continuity of the functions $s
\mapsto\mu([0,s] \times[0,1]^d)$, $u \mapsto\mu([0,1] \times[0,u]
\times[0,1]^{d-1}),\ldots, u \mapsto\mu([0,1]^d \times[0,u])$ on
$[0,1]$. This concludes the proof for the case $\xi_{i,n}^{(m)} \ge-K$.

Let us now consider the general case. Let $Z_{i,n}^+=\max(\xi
_{i,n}^{(m)},0)$, $Z_{i,n}^- = \max(-\xi_{i,n}^{(m)},0)$, $K^+ = \Ex
(Z_{0,n}^+)$ and $K^- = \Ex(Z_{0,n}^-)$. Furthermore, define $\xi
_{i,n}^{(m),+} = Z_{i,n}^+- K^+$ and $\xi_{i,n}^{(m),-} = Z_{i,n}^- -
K^-$. Then, using the fact that $K^+ - K^- = 0$, we can write
\[
\xi_{i,n}^{(m)} = Z_{i,n}^+ - Z_{i,n}^- =
Z_{i,n}^+- K^+ - \bigl( Z_{i,n}^-- K^- \bigr) =
\xi_{i,n}^{(m),+} - \xi_{i,n}^{(m),-}.
\]
Setting
\[
\tilde\B_n^{(m),\pm} (s,\mathbf u) = n^{-1/2} \sum
_{i=1}^{\lfloor ns
\rfloor} \xi_{i,n}^{(m),\pm} \bigl
\{ \1(\mathbf U_i \le\mathbf u) - C(\mathbf u) \bigr\}, \qquad(s,\mathbf u)
\in[0,1]^{d+1}, %
\]
we\vspace*{1pt} obtain that $\tilde\B_n^{(m)} = \tilde\B_n^{(m),+} - \tilde\B
_n^{(m),-}$. The case treated above immediately yields asymptotic
equicontinuity of $\tilde\B_n^{(m),+}$ and of $\tilde\B_n^{(m),-}$,
which implies asymptotic equicontinuity of $\tilde\B_{n}^{(m)}$.
\end{pf}

\begin{pf*}{Proof of Theorem~\ref{teomultCLT}}
Weak convergence of the finite-dimension\-al distributions is
established in Lemma~\ref{lemfidi}. Asymptotic tightness of $\tilde
\B_n$ is a consequence of the weak convergence of~$\tilde\B_n$
to $\B_C$ in $\ell^\infty([0,1]^d)$, which follows from Theorem~1 in
B\"ucher \cite{Buc14}. From Lemma~\ref{lemasymequi}, we have that,
for any $m \in\{1,\ldots,M\}$, $\tilde\B_n^{(m)}$ is asymptotically
uniformly $\llVert  \cdot \rrVert  _1$-equicontinuous in probability. Together with
the fact that $[0,1]^{d+1}$ is totally bounded for $\llVert  \cdot \rrVert  _1$ and
Lemma~\ref{lemfidi}, we have, for instance, from Theorem 2.1 in Kosorok \cite
{Kos08}, that, for any $m \in\{1,\ldots,M\}$, $\tilde\B_n^{(m)}
\leadsto\B_C^{(m)}$ in $\ell^\infty([0,1]^d)$, which implies
asymptotic tightness of $\tilde\B_n^{(m)}$. The proof is complete as
marginal asymptotic tightness implies joint asymptotic tightness.
\end{pf*}


\section{Proof of Theorem~\texorpdfstring{\protect\ref{teostuterep}}{3.4}}
\label{proofstuterep}

The proof of Theorem~\ref{teostuterep} is based on the extended
continuous mapping theorem (van der Vaart and Wellner \cite{vanWel96}, Theorem 1.11.1). The
intuition of the proof is as follows: the aim is to construct suitable
maps $g_n$ and $g$ such that $g_n$ continuously converges to $g$ (i.e.,
$g_n(\alpha_n)$ converges uniformly to $g(\alpha)$ for all sequences
$\alpha_n$ converging uniformly to $\alpha$) and such that we may
conclude that, as a process indexed\vspace*{1pt} by $s,t,\mathbf u$, $\Cb_n(s,t,\mathbf u)
\approx g_n\{\tilde\B_n(t, \mathbf u)- \tilde\B_n(s,\mathbf u) \}$
converges weakly to $g\{\tilde\B(t, \mathbf u)- \tilde\B(s,\mathbf u) \} =
\Cb(s,t,\mathbf u)$.

In the following, all the convergences are with respect to $n \to
\infty$. Let $\EE$ be the set of c.d.f.s on $[0,1]$ with no mass at
0, that is,
\begin{eqnarray*}
\EE &=& \bigl\{F\dvtx [0,1] \to[0,1]:   F \mbox{ is right-continuous and nondecreasing
with }
\\
&&\hspace*{173pt} F(0)=0 \mbox{ and } F(1) = 1 \bigr\},
\end{eqnarray*}
let
\begin{eqnarray*}
\EE_n^\star &=& \bigl\{F^\star\dvtx \Delta\times[0,1]
\to[0,1]:   u \mapsto \lambda_n(s,t)^{-1}
F^\star(s,t,u) \in\EE\mbox{ if } \lfloor ns \rfloor < \lfloor nt \rfloor
\\
&&\hspace*{153pt}\mbox{and } F^\star(s,t,\cdot) = 0\mbox{ if } \lfloor ns \rfloor =
\lfloor nt \rfloor \bigr\},
\end{eqnarray*}
where $\lambda_n(s,t) = (\lfloor nt \rfloor-\lfloor ns \rfloor) /
n$, and let $I_n$ be the sequence of maps defined, for any $F^\star\in
\EE_n^\star$ and any $(s,t,u) \in\Delta\times[0,1]$, by
\[
I_n\bigl(F^\star\bigr) (s,t,u) = \inf\bigl\{v \in[0,1]:
F^\star(s,t,v) \geq \lambda_n(s,t) u \bigr\}. %
\]
Furthermore, given a function $H^\star\in\ell^\infty(\Delta\times
[0,1]^d)$, for any $j \in\{1,\ldots,d\}$, we define
\[
H^\star_j(s,t,u) = H^\star(s,t,\mathbf
u_{\{j\}}), \qquad(s,t,u) \in \Delta\times[0,1], %
\]
where, for any $u \in[0,1]$, $\mathbf u_{\{j\}}$ is the vector of
$[0,1]^d$ whose components are all equal to 1 except the $j$th one
which is equal to $u$. Then, let
\[
\EE_{n,d}^\star= \bigl\{H^\star\dvtx \Delta
\times[0,1]^d \to[0,1]\dvtx  H_j^\star\in
\EE_n^\star\mbox{ for all } j \in\{1,\ldots,d\} \bigr\}
\]
and let $\Phi_n$ be the map from $\EE_{n,d}^\star$ to $\ell^\infty
(\Delta\times[0,1]^d)$ defined, for any $H^\star\in\EE_{n,d}^\star
$ and $(s,t,\mathbf u) \in\Delta\times[0,1]^d$, by
%
\begin{eqnarray}
\label{eqphin} \Phi_n\bigl(H^\star\bigr) (s,t,\mathbf u) =
H^\star\bigl\{ s,t,I_n\bigl(H_1^\star
\bigr) (s,t,u_1),\ldots,I_n\bigl(H_d^\star
\bigr) (s,t,u_d) \bigr\}.
\end{eqnarray}

Let additionally $U_n^\star\in\EE_n^\star$ be defined as $U_n^\star
(s,t,u) = \lambda_n(s,t) u$ for all $(s,t,u) \in\Delta\times[0,1]$,
and let $C_n^\star(s,t,\mathbf u) = \lambda_n(s,t) C(\mathbf u)$ for all
$(s,t,\mathbf u) \in\Delta\times[0,1]^d$. Clearly, we have that
$C_{n,1}^\star= \cdots= C^\star_{n,d} = U_n^\star$. Moreover, $\Phi
_n(C_n^\star)=C_n^\star$.

Also, let
\begin{eqnarray*}
\DD^\star&=& \bigl\{\alpha^\star\in\ell^\infty\bigl(
\Delta\times[0,1]^d\bigr):   \alpha^\star(s,t,\cdot) = 0
\mbox{ if } s = t, \mbox{ and}
\\
&&\hspace*{5pt} \alpha^\star(s,t,\mathbf u) = 0 \mbox{ if } s < t \mbox{ and if one of the
components of } \mathbf u \mbox{ is 0 or } \mathbf u = (1,\ldots,1) \bigr\},
\end{eqnarray*}
let $\DD_n^\star= \{\alpha^\star\in\DD^\star:   C_n^\star
+n^{-1/2} \alpha^\star\in\EE_{n,d}^\star\}$, and let $\DD_0^\star
= \DD^\star\cap \CC(\Delta\times[0,1]^d)$. Finally, for any
$\alpha_n^\star\in\DD_n^\star$ and any $(s,t,\mathbf u) \in\Delta
\times[0,1]^d$, let
%
\begin{eqnarray}
\label{eqgn} g_n\bigl(\alpha_n^\star\bigr)
(s,t,\mathbf u) = \sqrt n \bigl\{\Phi_n\bigl(C_n^\star
+n^{-1/2} \alpha_n^\star\bigr) (s,t,\mathbf u) -
\Phi_n\bigl(C_n^\star\bigr) (s,t,\mathbf u) \bigr\},
\end{eqnarray}
and, for any $\alpha^\star\in\DD_0^\star$ and any $(s,t,\mathbf u) \in
\Delta\times[0,1]^d$, let
\begin{eqnarray*}
g\bigl(\alpha^\star\bigr) (s,t,\mathbf u)= \alpha^\star(s,t,\mathbf u)
- \sum_{j=1}^d \dot C_j(\mathbf
u) \alpha^\star\bigl(s,t,\mathbf u^{(j)}\bigr).
\end{eqnarray*}
The following lemma is the main ingredient for the proof of Theorem
\ref{teostuterep}. Its proof is given subsequent to the proof of
Theorem \ref{teostuterep}.

%
\begin{lem}
\label{lemgn}
Suppose that $C$ satisfies Condition \ref{condpd}, and let $\alpha
^\star_n \to\alpha^\star$ with $\alpha^\star_n \in\DD_n^\star$
for every $n$ and $\alpha^\star\in\DD_0^\star$. Then, $g_n(\alpha
^\star_n) \to g(\alpha^\star) \in\ell^\infty(\Delta\times[0,1]^d)$.
\end{lem}

\begin{pf*}{Proof of Theorem~\ref{teostuterep}}
Under Condition~\ref{condweak}, we have that $\tilde\B_n \leadsto
\B_C$ in $\ell^\infty([0,1]^{d+1})$. Now, for any $(s,t,\mathbf u) \in
\Delta\times[0,1]^d$, define $\tilde\B_n^\Delta(s,t,\mathbf u) =
\tilde\B_n(t,\mathbf u) - \tilde\B_n(s,\mathbf u)$, $\B_C^\Delta(s,t,\mathbf
u) = \B_C(t,\mathbf u) - \B_C(s,\mathbf u)$, and
\[
\tilde H_n^\star(s,t,\mathbf u) = \frac{1}{n} \sum
_{i=\lfloor ns \rfloor
+1}^{\lfloor nt \rfloor} \1(\mathbf U_i \leq\mathbf
u). %
\]
Notice that $\tilde\B_n^\Delta= \sqrt{n} (\tilde H_n^\star-
C_n^\star)$ and that, by the continuous mapping theorem, $\tilde\B
_n^\Delta\leadsto\B_C^\Delta$ in $\ell^\infty(\Delta\times
[0,1]^d)$. Clearly, $\tilde\B_n^\Delta$, as a function of $\omega$,
takes its values in $\DD_n^\star$ and $\B_C^\Delta$ is Borel
measurable and separable by Condition~\ref{condweak}, and, as a
function of $\omega$, takes its values in $\DD_0^\star$. Now,
consider the map $h_n$ from $\DD_n^\star$ to $\{ \ell^\infty(\Delta
\times[0,1]^d) \}^2$, defined, for any $\alpha_n^\star\in\DD
_n^\star$ and any $(s,t,\mathbf u) \in\Delta\times[0,1]^d$, by
\[
h_n\bigl(\alpha_n^\star\bigr) (s,t,\mathbf u) =
\bigl(g_n\bigl(\alpha_n^\star\bigr) (s,t,\mathbf u),
g\bigl(\alpha_n^\star\bigr) (s,t,\mathbf u) \bigr). %
\]
Using Lemma~\ref{lemgn} and the fact that $g$ is linear and bounded,
we have from the extended continuous mapping theorem
(van der Vaart and Wellner \cite{vanWel96}, Theorem 1.11.1) that $h_n(\tilde\B_n^\Delta) \leadsto h(\B
_C^\Delta)$ in $\{ \ell^\infty(\Delta\times[0,1]^d) \}^2$, where,
for any $\alpha^\star\in\DD_0^\star$ and any $(s,t,\mathbf u) \in
\Delta\times[0,1]^d$,
\[
h\bigl(\alpha^\star\bigr) (s,t,\mathbf u) = \bigl(g\bigl(
\alpha^\star\bigr) (s,t,\mathbf u), g\bigl(\alpha^\star\bigr) (s,t,
\mathbf u) \bigr). %
\]
An application of the continuous mapping theorem immediately yields
that $g_n(\tilde\B_n^\Delta) - \tilde\Cb_n = g_n(\tilde\B
_n^\Delta) - g(\tilde\B_n^\Delta) \leadsto0$ in $\ell^\infty
(\Delta\times[0,1]^d)$, where $\tilde\Cb_n$ is defined in~(\ref
{eqtildeCn}). To complete the proof, it remains to show that
\[
A_n = \sup_{(s,t,\mathbf u) \in\Delta\times[0,1]^d} \bigl\llvert g_n
\bigl(\tilde \B_n^\Delta\bigr) (s,t,\mathbf u) -
\Cb_n(s,t,\mathbf u) \bigr\rrvert =\mathrm{o}_\Pr(1). %
\]
Note that it suffices to restrict the supremum over all pairs $(s,t)
\in\Delta$ such that $\lfloor ns \rfloor < \lfloor nt \rfloor$.
From the definition of $g_n$, we have that
\begin{eqnarray*}
&& g_n\bigl(\tilde\B_n^\Delta\bigr) (s,t,\mathbf u)
\\
&&\quad =
\sqrt{n} \bigl\{ \Phi _n\bigl(\tilde H_n^\star
\bigr) (s,t,\mathbf u) - \Phi_n\bigl(C_n^\star\bigr)
(s,t,\mathbf u) \bigr\}
\\
&&\quad = \frac{1}{\sqrt{n}} \sum_{i=\lfloor ns \rfloor+1}^{\lfloor nt
\rfloor} \bigl[
\1 \bigl\{ U_{i1} \leq I_n\bigl(\tilde
H_{n,1}^\star \bigr) (s,t,u_1),\ldots,
U_{id} \leq I_n\bigl(\tilde H_{n,d}^\star
\bigr) (s,t,u_d)\bigr\} - C(\mathbf u) \bigr].
\end{eqnarray*}
Now, let $\tilde H_{\lfloor ns \rfloor+1:\lfloor nt \rfloor}$ be the
empirical c.d.f. computed from the sample $\mathbf U_{\lfloor ns \rfloor
+1},\ldots,\mathbf U_{\lfloor nt \rfloor}$, 
and let $\tilde H_{\lfloor ns \rfloor+1:\lfloor nt \rfloor,1},\ldots,\tilde H_{\lfloor ns \rfloor+1:\lfloor nt \rfloor,d}$ be the
corresponding marginal c.d.f.s. Given $F \in\EE$, let $F^{-1}$ be its
generalized inverse defined by $F^{-1}(u) = \inf\{ v \in[0,1]:   F(v)
\geq u \}$. Then, let
\[
\tilde{\mathbf  H}_{\lfloor ns \rfloor+1:\lfloor nt \rfloor}^{-1}(\mathbf u) = \bigl( \tilde
H_{\lfloor ns \rfloor+1:\lfloor nt \rfloor
,1}^{-1}(u_1),\ldots,\tilde
H_{\lfloor ns \rfloor+1:\lfloor nt \rfloor
,d}^{-1}(u_d) \bigr), \qquad\mathbf u
\in[0,1]^d. %
\]
Using the fact that, for any $j \in\{1,\ldots,d\}$, $I_n(\tilde
H_{n,j}^\star)(s,t,u) = \tilde H_{\lfloor ns \rfloor+1:\lfloor nt
\rfloor,j}^{-1}(u)$ for all $(s,t,u) \in\Delta\times[0,1]$ such
that $\lfloor ns \rfloor < \lfloor nt \rfloor$, we obtain
\begin{eqnarray*}
g_n\bigl(\tilde\B_n^\Delta\bigr) (s,t,\mathbf u) &=&
\frac{1}{\sqrt{n}} \sum_{i=\lfloor ns \rfloor+1}^{\lfloor nt \rfloor} \bigl[ \1
\bigl\{ \mathbf U_i \leq \tilde{\mathbf  H}_{\lfloor ns \rfloor+1:\lfloor nt \rfloor
}^{-1}(\mathbf
u)\bigr\} - C(\mathbf u) \bigr]
\\
&=& \sqrt{n} \lambda_n(s,t) \bigl[ \tilde H_{\lfloor ns \rfloor
+1:\lfloor nt \rfloor} \bigl\{ \tilde{\mathbf  H}_{\lfloor ns \rfloor
+1:\lfloor nt \rfloor}^{-1}(\mathbf u) \bigr\} - C(\mathbf u) \bigr].
\end{eqnarray*}
Hence, we obtain that
\begin{eqnarray*}
A_n &=& \sup_{(s,t,\mathbf u) \in\Delta\times[0,1]^d} \sqrt{n} \lambda
_n(s,t) \bigl\llvert C_{\lfloor ns \rfloor + 1:\lfloor nt \rfloor}(\mathbf u) - \tilde
H_{\lfloor ns \rfloor+1:\lfloor nt \rfloor} \bigl\{ \tilde{\mathbf  H}_{\lfloor ns \rfloor+1:\lfloor nt \rfloor}^{-1}(\mathbf u)
\bigr\} \bigr\rrvert
\\
&=& n^{-1/2} \max_{1 \leq l < k \leq n} \sup_{\mathbf u \in[0,1]^d}
(k-l) \bigl\llvert C_{l + 1:  k}(\mathbf u) - \tilde H_{l+1:  k} \bigl\{ \tilde{\mathbf  H}_{l+1:  k}^{-1}(\mathbf u) \bigr\} \bigr\rrvert.
\end{eqnarray*}
Under Condition~\ref{condnoTies}, it can be verified, using
properties of generalized inverses, that
\[
\sup_{\mathbf u \in[0,1]^d} \bigl\llvert C_{l + 1:  k}(\mathbf u) - \tilde
H_{l+1:  k} \bigl\{ \tilde{\mathbf  H}_{l+1:  k}^{-1}(\mathbf u)
\bigr\} \bigr\rrvert \leq\frac{d}{k-l}, %
\]
which implies that $A_n \to0$ and completes the proof.
\end{pf*}


It remains to prove Lemma~\ref{lemgn}. For that purpose, another
lemma is needed.


%
\begin{lem} \label{leminv}
Let $\alpha^\star_n \to\alpha^\star$ with $\alpha^\star_n \in
\DD_n^\star$ for every $n$ and $\alpha^\star\in\DD_0^\star$.
Then, for any $j \in\{1,\ldots,d\}$,
\[
\sup_{(s,t,u) \in\Delta\times[0,1]} \bigl\llvert \sqrt{n} \lambda _n(s,t)
\bigl\{ I_n\bigl(U_n^\star+ n^{-1/2}
\alpha_{n,j}^\star\bigr) (s,t,u) - u \bigr\} +
\alpha_j^\star(s,t,u) \bigr\rrvert \to0. %
\]
\end{lem}

\begin{pf}
The assertion is trivial for $u = 0$ because $\alpha^\star\in\DD
_0^\star$ and $U_n^\star+n^{-1/2}\alpha_{n,j}^\star\in\EE_n^\star$.

Clearly, for any $s \in[0,1]$, $ns \geq\lfloor ns \rfloor$, that is,
$s \geq\lambda_n(0,s)$. Furthermore, under the constraint $s \leq t$,
$\lfloor nt \rfloor = \lfloor ns \rfloor$ is equivalent to $0 \leq
t-\lambda_n(0,s) < 1/n$, which can be written as $0 \leq t-s +
s-\lambda_n(0,s) < 1/n$, which means that there exists $h_n \downarrow
0$ such that $t-s < h_n$. Then, we have
\begin{eqnarray*}
&& \sup_{\lfloor nt \rfloor = \lfloor ns \rfloor, u \in[0,1]} \bigl\llvert \lambda _n(s,t) \sqrt{n}
\bigl\{ I_n\bigl(U_n^\star+ n^{-1/2}
\alpha_{n,j}^\star \bigr) (s,t,u) - u \bigr\} +
\alpha_j^\star(s,t,u) \bigr\rrvert
\\
&&\quad \leq\sup_{t-s <
h_n, u \in[0,1]} \bigl\llvert \alpha_j^\star(s,t,u)
\bigr\rrvert \to0
\end{eqnarray*}
by uniform continuity of $\alpha_j^\star$ on $\Delta\times[0,1]$.

Hence, it remains to consider the case $\lfloor ns \rfloor<\lfloor nt
\rfloor$ and $u \in(0,1]$. Given $F \in\EE$, let $F^{-1}$ be its
generalized inverse\vspace*{1pt} defined by $F^{-1}(u) = \inf\{ v \in[0,1]:   F(v)
\geq u \}$. Then, notice that, for any $\lfloor ns \rfloor < \lfloor
nt \rfloor$ and $u \in[0,1]$, $I_n(U_n^\star+ n^{-1/2} \alpha
_{n,j}^\star)(s,t,u) = F_{s,t,n}^{-1}(u)$, where $F_{s,t,n} = \lambda
_n(s,t)^{-1} (U_n^\star+ n^{-1/2} \alpha_{n,j}^\star)(s,t,\cdot)
\in\EE$. It follows\vspace*{1pt} that, for any $\lfloor ns \rfloor < \lfloor nt
\rfloor$ and $u \in(0,1]$, $\xi_n(s,t,u) = I_n(U_n^\star+ n^{-1/2}
\alpha_{n,j}^\star)(s,t,u) >0$, and therefore that $\varepsilon
_n(s,t,u) = n^{-1} \wedge\xi_n(s,t,u)>0$. Also, for any $F \in\EE$,
it can be verified that $F \{ F^{-1}(u) - \eta\} \leq u \leq F \circ
F^{-1}(u)$ for all $u \in(0,1]$ and all $\eta> 0$ such that
$F^{-1}(u) - \eta\ge0$.
Hence, for any $\lfloor ns \rfloor < \lfloor nt \rfloor$ and $u \in(0,1]$,
\begin{eqnarray*}
\bigl(U_n^\star+n^{-1/2}\alpha_{n,j}^\star
\bigr) \bigl\{ s, t,\xi_n(s,t,u) - \varepsilon_n(s,t,u)
\bigr\} &\le&\lambda_n(s,t)u
\\
&\le&\bigl(U_n^\star
+n^{-1/2}\alpha_{n,j}^\star\bigr) \bigl\{ s,t,
\xi_n(s,t,u) \bigr\}, %
\end{eqnarray*}
that is
%
\begin{eqnarray}
\label{eqestxi}
&& - n^{-1/2} \alpha_{n,j}^\star\bigl\{
s,t, \xi_n(s,t,u) \bigr\}\nonumber
\\
&&\quad \le\lambda _n(s,t)\bigl\{
\xi_n(s,t,u) - u \bigr\}
\\
&&\quad \le\lambda_n(s,t) \varepsilon_n(s,t,u) -
n^{-1/2}\alpha_{n,j}^\star \bigl\{ s,t,
\xi_n(s,t,u) - \varepsilon_n(s,t,u) \bigr\},\nonumber
\end{eqnarray}
which in turn implies that
%
\begin{equation}
\label{intermediate} \sup_{\lfloor ns \rfloor < \lfloor nt \rfloor, u \in(0,1]} \bigl\llvert \lambda_n(s,t)
\bigl\{\xi_n(s,t,u) - u\bigr\} \bigr\rrvert \to0
\end{equation}
since, by uniform convergence of $\alpha_n^\star$ to $\alpha^\star$
and the fact that $\alpha^\star\in\DD_0^\star$, the quantity $\sup_{(s,t,u) \in\Delta\times[0,1]} \llvert   \alpha_{n,j}^\star(s,t,u)\rrvert  $ is
bounded. From~(\ref{eqestxi}), exploiting the fact that $\varepsilon
_n(s,t,u) \le n^{-1}$, we then obtain that
\[
\sup_{\lfloor ns \rfloor < \lfloor nt \rfloor, u \in(0,1]} \bigl\llvert \sqrt{n} \lambda_n(s,t)
\bigl\{\xi_n(s,t,u) - u \bigr\} + \alpha_j^\star
(s,t,u) \bigr\rrvert 
\le A_n + B_n +
n^{-1/2},
\]
where
\[
A_n = \sup_{\lfloor ns \rfloor < \lfloor nt \rfloor, u \in(0,1]} \bigl\llvert
\alpha_n^\star\bigl\{s,t, \xi_n(s,t,u)\bigr\} -
\alpha_j^\star (s,t,u) \bigr\rrvert, %
\]
and
\[
B_{n} = \sup_{\lfloor ns \rfloor < \lfloor nt \rfloor, u \in(0,1]} \bigl\llvert
\alpha_{n,j}^\star\bigl\{s,t, \xi_n(s,t,u) -
\varepsilon _n(s,t,u)\bigr\} - \alpha_j^\star(s,t,u)
\bigr\rrvert. %
\]
%
For $B_{n}$, we write $B_{n} \leq B_{n,1} + B_{n,2}$, where
\begin{eqnarray*}
B_{n,1} &=& \mathop{\sup_{\lfloor ns \rfloor < \lfloor nt \rfloor}}_{u \in
(0,1]} \bigl\llvert
\alpha_{n,j}^\star\bigl\{s,t, \xi_n(s,t,u) -
\varepsilon _n(s,t,u)\bigr\} - \alpha_j^\star
\bigl\{s,t, \xi_n(s,t,u) - \varepsilon _n(s,t,u)\bigr\}
\bigr\rrvert
\\
&\leq&\sup_{(s,t,u) \in\Delta\times[0,1]} \bigl\llvert \alpha_{n,j}^\star(s,t,u)
- \alpha_j^\star(s,t,u) \bigr\rrvert \to0,
\end{eqnarray*}
and
\[
B_{n,2} = \sup_{(s, t, u) \in\Delta\times[0,1]} \bigl\llvert \alpha
_j^\star\bigl\{s,t, \xi_n(s,t,u) -
\varepsilon_n(s,t,u)\bigr\} - \alpha _j^\star(s,t,u)
\bigr\rrvert. %
\]
It remains to show that $B_{n,2} \to0$. Let $\varepsilon> 0$. Since
$\alpha^\star\in\DD_0^\star$, there exists $\delta> 0$ such that
$\sup_{t - s < \delta, u \in[0,1]} \llvert   \alpha_j^\star(s,t,u) \rrvert   \leq
\varepsilon$. We have $B_{n,2} = \max\{ B_{n,3}, B_{n,4} \}$, where
\[
B_{n,3} = \sup_{t - s < \delta, u \in[0,1]} \bigl\llvert
\alpha_j^\star \bigl\{s,t, \xi_n(s,t,u) -
\varepsilon_n(s,t,u)\bigr\} - \alpha_j^\star
(s,t,u) \bigr\rrvert \leq2\varepsilon, %
\]
and
\[
B_{n,4} = \sup_{t - s \geq\delta, u \in[0,1]} \bigl\llvert \alpha
_j^\star\bigl\{s,t, \xi_n(s,t,u) -
\varepsilon_n(s,t,u)\bigr\} - \alpha _j^\star(s,t,u)
\bigr\rrvert. %
\]
Now, it is easy to verify that $t-s \leq\lambda_n(s,t) + 1/n$, so
that, for $n$ sufficiently large, $t-s \geq\delta$ implies that
$\lambda_n(s,t) \geq\delta/2$. Then, from~(\ref{intermediate}) and
the fact that $\xi_n(\cdot,\cdot,0)=0$, we immediately have that,
for $n$ sufficiently large,
\[
a_n = \mathop{\sup_{t - s \geq\delta}}_{ u \in[0,1]} \bigl\llvert
\xi_n(s,t,u) - u \bigr\rrvert \leq\mathop{\sup_{t - s \geq\delta}}_{u \in[0,1]} \bigl
\llvert \lambda_n(s,t) \bigl\{ \xi_n(s,t,u) - u\bigr\}
\bigr\rrvert \times \sup_{t - s \geq\delta} \lambda _n(s,t)^{-1}
\to0. %
\]
Hence, we can write
\[
B_{n,4} \leq\mathop{\sup_{t - s \geq\delta, u, u' \in[0,1] }}_{\llvert  u' - u \rrvert
\leq a_n + n^{-1}} \bigl\llvert
\alpha_j^\star\bigl(s,t,u'\bigr) -
\alpha_j^\star (s,t,u) \bigr\rrvert \to0 %
\]
since\vspace*{1pt} $\alpha_j^\star$ is uniformly continuous on $\Delta\times
[0,1]$. Proceeding as for $B_{n}$, it can be verified that $A_n \to0$,
which completes the proof.
\end{pf}


\begin{pf*}{Proof of Lemma \ref{lemgn}}
Starting from the definitions of $g_n$ and $\Phi_n$ given in~(\ref
{eqgn}) and~(\ref{eqphin}), respectively, we have the decomposition
\[
g_n\bigl(\alpha_n^\star\bigr) (s,t,\mathbf u) =
A_{n,1}(s,t,\mathbf u) + A_{n,2}(s,t,\mathbf u), %
\]
where
\[
A_{n,1}(s,t,\mathbf u) = \alpha_n ^\star\bigl\{ s,t,
I_n\bigl(U_n^\star +n^{-1/2}
\alpha_{n,1}^\star\bigr) (s,t,u_1),\ldots,
I_n\bigl(U_n^\star +n^{-1/2}
\alpha_{n,d}^\star\bigr) (s,t,u_d) \bigr\},
\]
and
\begin{eqnarray*}
&& A_{n,2}(s,t,\mathbf u)
\\
&& = \sqrt{n} \lambda_n(s,t) \bigl[ C\bigl\{ I_n
\bigl(U_n^\star+n^{-1/2}\alpha_{n,1}^\star
\bigr) (s,t,u_1),\ldots, I_n\bigl(U_n^\star+n^{-1/2}
\alpha_{n,d}^\star\bigr) (s,t,u_d) \bigr\} - C(\mathbf
u) \bigr].
\end{eqnarray*}

We begin the proof by showing that $\sup_{(s,t,\mathbf u) \in\Delta
\times[0,1]^d} \llvert   A_{n,1}(s,t,\mathbf u) - \alpha^\star(s,t,\mathbf u) \rrvert   \to
0$. Let $\varepsilon>0$. Using the fact that $\alpha^\star\in\DD
_0^\star$, there exists $\delta> 0$ such that $\llvert  \alpha^\star
(s,t,\mathbf u)\rrvert   \le\varepsilon$ for all $t-s< \delta$ and $\mathbf u \in
[0,1]^d$. Then, we write
\[
\sup_{(s,t,\mathbf u) \in\Delta\times[0,1]^d} \bigl\llvert A_{n,1}(s,t,\mathbf u) -
\alpha^\star(s,t,\mathbf u) \bigr\rrvert \leq B_{n,1} +
B_{n,2} + B_{n,3},
\]
where
\begin{eqnarray*}
B_{n,1} &=&\sup_{(s,t,\mathbf u) \in\Delta\times[0,1]^d} \bigl\llvert
A_{n,1}(s,t,\mathbf u) - \alpha^\star\bigl\{ s,t, I_n
\bigl(U_n^\star +n^{-1/2}\alpha_{n,1}^\star
\bigr) (s,t,u_1),\ldots,
\\
&&{}\hspace*{172pt} I_n\bigl(U_n^\star
+n^{-1/2}\alpha_{n,d}^\star\bigr) (s,t,u_d)
\bigr\} \bigr\rrvert
\\
&\leq&\sup_{(s,t,\mathbf u) \in\Delta\times[0,1]^d} \bigl\llvert \alpha_n^\star
(s,t,\mathbf u) - \alpha^\star(s,t,\mathbf u) \bigr\rrvert \leq\varepsilon,
\end{eqnarray*}
for sufficiently large $n$, where
\begin{eqnarray*}
B_{n,2} &=& \mathop{\sup_{t - s < \delta}}_{\mathbf u \in[0,1]^d} \bigl\llvert \alpha
^\star\bigl\{ s,t, I_n\bigl(U_n^\star+n^{-1/2}
\alpha_{n,1}^\star\bigr) (s,t,u_1),\ldots,
I_n\bigl(U_n^\star+n^{-1/2}
\alpha_{n,d}^\star\bigr) (s,t,u_d) \bigr\}
\\
&&\hspace*{270pt}{} -
\alpha^\star(s,t,\mathbf u) \bigr\rrvert,
\end{eqnarray*}
and
\begin{eqnarray*}
B_{n,3} &=& \mathop{\sup_{t - s \geq\delta}}_{\mathbf u \in[0,1]^d}\bigl\llvert \alpha
^\star\bigl\{ s,t, I_n\bigl(U_n^\star+n^{-1/2}
\alpha_{n,1}^\star\bigr) (s,t,u_1),\ldots,
I_n\bigl(U_n^\star+n^{-1/2}
\alpha_{n,d}^\star\bigr) (s,t,u_d) \bigr\}
\\
&&\hspace*{270pt}{} -
\alpha^\star(s,t,\mathbf u) \bigr\rrvert.
\end{eqnarray*}
For $B_{n,2}$, using the triangle inequality, we have that
\[
B_{n,2} \leq2 \sup_{t-s <\delta, \mathbf u \in[0,1]^d} \bigl\llvert \alpha
^\star(s,t,\mathbf u) \bigr\rrvert \leq2 \varepsilon.
\]
%
For $B_{n,3}$, we use the fact that Lemma \ref{leminv} implies that,
for any $j \in\{1,\ldots,d\}$,
%
\begin{equation}
\label{eqanj} a_{n,j} = \sup_{t-s \ge\delta, u \in[0,1]} \bigl\llvert
I_n\bigl(U_n^\star+ n^{-1/2}
\alpha_{n,j}^\star\bigr) (s,t,u) - u \bigr\rrvert \to0,
\end{equation}
%
and the fact that
\[
B_{n,3} \leq\sup_{t - s \geq\delta, \llvert  u_1-v_1\rrvert   \leq a_{n,1},\ldots,
\llvert  u_d-v_d\rrvert   \leq a_{n,d} } \bigl\llvert
\alpha^\star(s,t, \mathbf u) - \alpha ^\star(s,t,\mathbf v) \bigr\rrvert. %
\]
By uniform continuity of $\alpha^\star$, for sufficiently large $n$,
we obtain that $B_{n,3} \leq\varepsilon$. Hence, we have shown that,
for sufficiently large $n$, $\sup_{(s,t,\mathbf u)\in\Delta\times
[0,1]^d} \llvert   A_{n,1}(s,t,\mathbf u) - \alpha^\star(s,t,\mathbf u) \rrvert   \leq4
\varepsilon$, and therefore that $\sup_{(s,t,\mathbf u)\in\Delta\times
[0,1]^d} \llvert   A_{n,1}(s,t,\mathbf u) - \alpha^\star(s,t,\mathbf u) \rrvert   \to0$.

Let us now deal with $A_{n,2}$. Fix $n \ge1$ and $s < t$ such that
$\lfloor ns \rfloor < \lfloor nt \rfloor$. For any $\mathbf u \in
[0,1]^d$, $j \in\{1,\ldots,d\}$ and $r \in[0,1]$, let $\bar u_j(r) =
u_j + r \{ I_n(U_n^\star+n^{-1/2}\alpha_{n,j}^\star)(s,t,u_j) - u_j
\}$ and define $\bar{\mathbf u} (r) = ( \bar u_1(r),\ldots,\bar u_d(r) )$.
Now, fix $\mathbf u \in(0,1)^d$ and let $f$ be the function defined by
\[
f(r) = C_n^\star\bigl\{ s,t,\bar{\mathbf u} (r) \bigr\} =
\lambda_n(s,t) C \bigl\{ \bar {\mathbf u} (r) \bigr\}. %
\]
Obviously, we have that $0 < \bar u_j(r) < 1$ for all $r \in(0,1)$ and
$j \in\{1,\ldots,d\}$. Therefore, the function $f$ is continuous on
$[0,1]$, and, by Condition~\ref{condpd}, is differentiable on
$(0,1)$. Hence, by the mean value theorem, there exists $r^* \in(0,1)$
such that $f(1) - f(0) = f'(r^*)$, which implies that
%
\begin{equation}
\label{eqMVT} A_{n,2}(s,t,\mathbf u) = \sum_{j=1}^d
\dot C_j\bigl\{ \bar{\mathbf u} \bigl(r^*\bigr) \bigr\}
\lambda_n(s,t) \sqrt{n} \bigl\{ I_n\bigl(U_n^\star+
n^{-1/2} \alpha _{n,j}^\star\bigr)
(s,t,u_j) - u_j \bigr\}.
\end{equation}

The previous equality remains clearly valid when $\lfloor ns \rfloor =
\lfloor nt \rfloor$. Let us now verify that it also holds when
$\lfloor ns \rfloor < \lfloor nt \rfloor$ and $\mathbf u$ is on the
boundary of $[0,1]^d$. When $u_j=0$ for some $j \in\{1,\ldots,d\}$,
$I_n(U_n^\star+ n^{-1/2} \alpha_{n,j}^\star)(\cdot,\cdot,u_j) =
0$, which implies that $\bar u_j(r) = 0$ for all $r \in[0,1]$. It then
immediately follows that the left-hand side of~(\ref{eqMVT}) is zero and
that the $j$th term in the sum on the right is zero. The $d-1$
remaining terms in the sum on the right of~(\ref{eqMVT}) are actually
also zero because, for any $k \in\{1,\ldots,d\}$, $k \neq j$, $\dot
C_k(\mathbf v) = 0$ for all $\mathbf v \in[0,1]^d$ such that $v_k = 0$.
Hence,~(\ref{eqMVT}) remains true whenever $u_j = 0$ for some $j \in
\{1,\ldots,d\}$.

Let us now assume that $\lfloor ns \rfloor < \lfloor nt \rfloor$ and
that $u_j=1$ for some $j \in\{1,\ldots,d\}$. Two cases can be
distinguished according to whether $I_n(U_n^\star+ n^{-1/2} \alpha
_{n,j}^\star)(s,t,1)=1$ or $I_n(U_n^\star+ n^{-1/2} \alpha
_{n,j}^\star)(s,t,1) < 1$. In the\vspace*{1pt} later case, $0 < \bar u_j(r) < 1$.
In the former case, we obtain that $\bar u_j(r) = 1$ for all $r \in
[0,1]$ and that the $j$th term in the sum on the right of~(\ref
{eqMVT}) is zero so that neither the left nor the right-hand side of~(\ref
{eqMVT}) depend on $u_j$ anymore. It follows that, when some
components of $\mathbf u$ are one, the previous equality can be recovered
by an application of the mean value theorem similar to the one carried
out above.

Now, we write
%
\begin{equation}
\label{eqAn2} A_{n,2}(s,t,\mathbf u) = \sum_{j=1}^d
\dot C_j(\mathbf u) \lambda_n(s,t) \sqrt{n} \bigl\{
I_n\bigl(U_n^\star+ n^{-1/2}
\alpha_{n,j}^\star\bigr) (s,t,u_j) -
u_j \bigr\} + r_n(s,t,\mathbf u),\quad
\end{equation}
where $r_n(s,t,\mathbf u) = \sum_{j=1}^d r_{n,j}(s,t,\mathbf u)$ and, for any
$j \in\{1,\ldots,d\}$,
\[
r_{n,j}(s,t,\mathbf u) = \bigl[ \dot C_j\bigl\{\bar{\mathbf u}
\bigl(r^*\bigr)\bigr\} - \dot C_j(\mathbf u) \bigr] \lambda_n(s,t)
\sqrt{n} \bigl\{ I_n\bigl(U_n^\star+
n^{-1/2} \alpha _{n,j}^\star\bigr)
(s,t,u_j) - u_j \bigr\}. %
\]
By Lemma~\ref{leminv} and from the fact that $0 \leq\dot C_j \leq1$
for all $j \in\{1,\ldots,d\}$, the dominating term in
decomposition~(\ref{eqAn2}) converges to
\[
- \sum_{j=1}^d \dot C_j(\mathbf
u) \alpha^\star\bigl(s,t,\mathbf u^{(j)}\bigr) %
\]
uniformly in $(s, t, \mathbf u) \in\Delta\times[0,1]^d$. It therefore
remains to show that
\[
\sup_{(s, t, \mathbf u) \in\Delta\times[0,1]^d} \bigl\llvert r_n(s,t,\mathbf u) \bigr
\rrvert \to0. %
\]

Let us first show that $\sup_{(s, t, \mathbf u) \in\Delta\times[0,1]^d}
\llvert   r_{n,1}(s,t,\mathbf u) \rrvert   \to0$. We have that
\[
\sup_{(s, t, \mathbf u) \in\Delta\times[0,1]^d} \bigl\llvert r_{n,1}(s,t,\mathbf u) \bigr
\rrvert \leq B_{n,4} + B_{n,5}, %
\]
where
\begin{eqnarray*}
B_{n,4} &=& \sup_{(s, t, \mathbf u) \in\Delta\times[0,1]^d} \bigl\llvert \dot
C_1\bigl\{ \bar{\mathbf u} \bigl(r^*\bigr)\bigr\} - \dot C_1(
\mathbf u) \bigr\rrvert
\\
&&{}\times \sup_{(s, t, \mathbf
u) \in\Delta\times[0,1]^d} \bigl\llvert \lambda_n(s,t)
\sqrt{n} \bigl\{ I_n\bigl(U_n^\star+
n^{-1/2} \alpha_{n,1}^\star\bigr) (s,t,u_1)
- u_1 \bigr\} + \alpha_1^\star(s,t,u_1)
\bigr\rrvert,
\end{eqnarray*}
and
\[
B_{n,5} = \sup_{(s, t, \mathbf u) \in\Delta\times[0,1]^d} \bigl\llvert \bigl[ \dot
C_1\bigl\{\bar{\mathbf u} \bigl(r^*\bigr)\bigr\} - \dot C_1(
\mathbf u) \bigr] \alpha_1^\star (s,t,u_1) \bigr
\rrvert. %
\]
From the fact that $0 \leq\dot C_1 \leq1$ and Lemma~\ref{leminv},
we immediately obtain that $B_{n,4} \to0$. It remains to show that
$B_{n,5} \to0$. To this end, let $\varepsilon>0$. Since $\alpha
^\star\in\DD_0^\star$, there exists $\delta>0 $ such that $\llvert  \alpha
_1^\star(s,t,u)\rrvert   \le\varepsilon$ for all $t-s < \delta$ and all $u
\in[0,1]$. Then, $B_{n,5} \leq B_{n,6} + B_{n,7}$, where
\[
B_{n,6} = \sup_{(s, t, \mathbf u) \in\Delta\times[0,1]^d} \bigl\llvert \dot
C_1\bigl\{ \bar{\mathbf u} \bigl(r^*\bigr)\bigr\} - \dot C_1(
\mathbf u) \bigr\rrvert \times \sup_{t - s < \delta, u \in[0,1]} \bigl\llvert
\alpha_1^\star(s,t,u) \bigr\rrvert \leq2 \varepsilon,
\]
and
\[
B_{n,7} = \sup_{t - s \geq\delta, \mathbf u \in[0,1]^d} \bigl\llvert \bigl[ \dot
C_1\bigl\{\bar{\mathbf u} \bigl(r^*\bigr)\bigr\} - \dot C_1(
\mathbf u) \bigr] \alpha_1^\star(s,t,u_1) \bigr
\rrvert. %
\]
For $B_{n,7}$, we use the fact that, since $\alpha^\star\in\DD
_0^\star$, there exists $0 < \kappa< 1/2$ such that
\[
\sup_{t - s \geq\delta, u \in[0,\kappa) \cup(1-\kappa,1]} \bigl\llvert \alpha_1^\star(s,t,u)
\bigr\rrvert \leq\varepsilon. %
\]
Then, we write $B_{n,7} \leq B_{n,8} + B_{n,9}$, where
\[
B_{n,8} = \sup_{(s,t,\mathbf u) \in\Delta\times[0,1]^d} \bigl\llvert \dot
C_1\bigl\{ \bar{\mathbf u} \bigl(r^*\bigr)\bigr\} - \dot C_1(
\mathbf u) \bigr\rrvert \times \mathop{\sup_{t - s \geq
\delta, \mathbf u \in[0,1]^d }}_{ u_1 \in[0,\kappa) \cup(1-\kappa,1]} \bigl\llvert
\alpha_1^\star(s,t,u_1) \bigr\rrvert \leq2
\varepsilon, %
\]
and
\[
B_{n,9} = \sup_{t - s \geq\delta, \mathbf u \in[0,1]^d, u_1 \in[\kappa, 1- \kappa]} \bigl\llvert \dot
C_1\bigl\{\bar{\mathbf u} \bigl(r^*\bigr)\bigr\} - \dot C_1(
\mathbf u) \bigr\rrvert \times\sup_{(s,t,u) \in\Delta\times[0,1]} \bigl\llvert
\alpha_1^\star (s,t,u) \bigr\rrvert. %
\]
From~(\ref{eqanj}), we obtain that
\[
B_{n,9} \leq\mathop{\sup_{\mathbf u, \mathbf v \in[0,1]^d, u_1,v_1 \in[\kappa/2,
1- \kappa/2] }}_{\llvert  u_1 - v_1\rrvert   \leq a_{n,1},\ldots, \llvert  u_d - v_d\rrvert   \leq
a_{n,d}} \bigl\llvert \dot
C_1(\mathbf u) - \dot C_1(\mathbf v) \bigr\rrvert \times \sup
_{(s,t,u)
\in\Delta\times[0,1]} \bigl\llvert \alpha_1^\star(s,t,u)
\bigr\rrvert. %
\]
Since $\dot C_1$ is uniformly continuous on $[\kappa/2, 1- \kappa/2]
\times[0,1]^{d-1}$ according to Condition~\ref{condpd}, and since
$\sup_{(s,t,u) \in\Delta\times[0,1]} \llvert   \alpha_1^\star
(s,t,u) \rrvert  $ is bounded, we have that $B_{n,9} \to0$, which
implies that, for $n$ sufficiently large, $B_{n,9} \leq\varepsilon$.
It follows that, for $n$ sufficiently large, $B_{n,5} \leq5
\varepsilon$, which implies that $\sup_{(s,t,\mathbf u) \in\Delta\times
[0,1]^d} \llvert   r_{n,1}(s,t,\mathbf u) \rrvert   \to0$. One can proceed similarly for
$r_{n,j}$, $j \in\{2,\ldots,d\}$. Hence, $\sup_{s \leq t, \mathbf u \in
[0,1]^d} \llvert   r_n(s,t,\mathbf u) \rrvert   \to0$.
\end{pf*}
\end{appendix}

\section*{Acknowledgements}
The authors would like to thank two anonymous referees and an Associate
Editor for their constructive comments on an earlier version of this
manuscript, as well as Johan Segers for fruitful discussions as usual.

This work has been supported in parts by the Collaborative Research
Center ``Statistical modeling of nonlinear dynamic processes'' (SFB
823) of the German Research Foundation (DFG) and by the IAP research
network Grant P7/06 of the Belgian government (Belgian Science Policy),
which is gratefully acknowledged.

\begin{supplement}
\stitle{Supplement to ``A dependent multiplier bootstrap for the sequential empirical
copula process under strong mixing''}
\slink[doi]{10.3150/14-BEJ682SUPP} 
\sdatatype{.pdf}
\sfilename{BEJ682\_supp.pdf}
\sdescription{Additional proofs and simulation results can be found in (B\"ucher and Kojadinovic \cite{BucKoj14}).}
\end{supplement}

%

\printhistory
\end{document}